\documentclass{statsoc}
\usepackage{mathptm}   
\usepackage{mathbbol}
\usepackage{mathrsfs}
\usepackage{comb, amsmath, latexsym, amssymb, amscd}
\usepackage{amsbsy}
\usepackage{natbib}
\usepackage{longtable}
\usepackage{graphicx}
\usepackage{dcpic}
\usepackage{pictex}

\def\widebar{\overline}

\title[An invitation to quantum tomography]{An invitation to quantum tomography}
\author{L.M. Artiles}
\address{Eurandom, P.O. Box 513, 5600 MB Eindhoven, 
The Netherlands, 
artiles@eurandom.tue.nl, http://euridice.tue.nl/$\sim$lartiles/} 
\author{R. Gill}
\address{Mathematical Institute, University of Utrecht, 
Box 80010, 3508 TA Utrecht, The Netherlands, gill@math.uu.nl, http://www.math.uu.nl/people/gill}
\author[L.M. Artiles {\it et al.}]{M.I. Gu\c{t}\u{a}}
\address{Eurandom, P.O. Box 513, 5600 MB Eindhoven, The Netherlands, 
guta@eurandom.tue.nl, http://euridice.tue.nl/$\sim$mguta/} 

\begin{document}
\DeclareGraphicsExtensions{.eps}
    
\maketitle


\begin{abstract}
The quantum state of a light beam can be represented as an infinite 
dimensional density matrix or equivalently as a density on the plane called the Wigner function. We describe quantum tomography as an inverse statistical problem in which the state is the unknown parameter and the data is given by results of measurements performed on identical quantum systems. 
We present consistency results for Pattern Function Projection Estimators as well 
as for Sieve Maximum Likelihood Estimators for both the density matrix of the quantum state and its Wigner function. Finally we illustrate via simulated data the performance of the estimators. An EM algorithm is proposed for practical implementation. 
There remain many open problems, e.g.~rates of convergence, adaptation, studying other estimators, etc., and a main purpose of the paper is to bring these to the attention 
of the statistical community.


\end{abstract}

\section{Introduction}
\label{sec.introduction}

It took more than eighty years from its discovery till it was possible 
to experimentally determine and visualize the most fundamental object in quantum mechanics, the wave function. The forward route from quantum state to probability distribution of measurement results has been the basic stuff of quantum mechanics textbooks for decennia. 
That the corresponding mathematical inverse problem had a solution, provided 
(speaking metaphorically) that the quantum state has been probed from a sufficiently rich set of directions, had also been known for many years. However it was only with  \cite{Smithey}, that it became feasible to actually carry out the corresponding measurements on one particular quantum system---in that case, the state of one mode of electromagnetic radiation (a pulse of laser light at a given frequency).  Experimentalists have used the technique to establish that they have succeeded in creating non-classical forms of laser 
light such as squeezed light and Schr\"odinger cats. The experimental technique we are referring to here is called quantum homodyne tomography: the word homodyne referring to a comparison between the light being measured with a reference light beam at the same frequency. We will explain the word tomography in a moment.

The quantum state can be represented mathematically in many different but equivalent ways, all of them linear transformations on one another. One favorite is as the Wigner function $W$: a real function of two variables, integrating to plus one over the whole plane, but not necessarily nonnegative. It can be thought of as a ``generalized joint probability density'' of the electric and magnetic fields, $q$ and $p$. However one cannot measure both fields at the same time and in quantum mechanics it makes no sense to talk about the values of both electric and magnetic fields simultaneously. It does, however, make sense to talk about the value of any linear combination of the two fields, say $\cos (\phi) q + \sin (\phi) p$. 
And one way to think about the statistical problem is as follows: 
the unknown parameter is a joint probability density $W$ of two variables $Q$ and $P$. The data consists of independent samples from the distribution of 
$(X,\Phi)=(\cos (\Phi) Q + \sin (\Phi) P,\Phi)$, where $\Phi$ is chosen independently of 
$(Q,P)$, and uniformly in the interval $[0,\pi]$. Write down the mathematical model 
expressing the joint density of $(X,\Phi)$ in terms of that of $(Q,P)$. Now just allow that 
latter joint density, $W$, to take negative as well as positive values (subject to certain 
restrictions which we will mention later). And that is the statistical problem of this paper.

This is indeed a classical tomography problem: we take observations from all possible 
one-dimensional projections of a two-dimensional density. The non-classical 
feature is that though all these one-dimensional projections are indeed bona-fide probability densities, the underlying two-dimensional ``joint density'' need not itself be a bona-fide joint probability density, but can have small patches of ``negative probability''.

Though the parameter to be estimated may look strange from some points of view, it is mathematically very nice from others. For instance, one can also represent it by a matrix of (a kind of) Fourier coefficients: one speaks then of the ``density matrix'' $\rho$. This is an infinite dimensional matrix of complex numbers, but it is a positive and selfadjoint matrix with trace one. The diagonal elements are real numbers summing up to one, and forming the probability distribution of the number of photons found in the light beam (if one could do that measurement). 
Conversely, any such matrix $\rho$ corresponds to a physically possible Wigner function $W$, so we have here a concise mathematical characterization of precisely which ``generalized joint probability densities'' can occur.

The initial reconstructions were done by borrowing analytic techniques from classical 
tomography -- the data was binned and smoothed, the inverse Radon transform carried out, followed by some Fourier transformations. At each of a number of steps, there are numerical discretization and truncation errors. The histogram of the data will not lie in the range of the forward transformation (from quantum state to density of the data). Thus the result of blindly applying an inverse will not be a bona-fide Wigner function or density matrix. Moreover the various numerical approximations all involve arbitrary choices of smoothing, binning or truncation parameters. Consequently the final picture can look just how the experimenter would like it to look and there is no way to statistically evaluate the reliability of the result. On the other hand the various numerical approximations tend to destroy the interesting ``quantum'' features the experimenter is looking for, so this method lost in popularity after the initial enthousiasm.

So far there has been little attention paid to this problem by statisticians, although 
on the one hand it is an important statistical problem coming up in 
modern physics, and on the other hand it is ``just'' a classical nonparametric statistical 
inverse problem. The unknown parameter is some object $\rho$, or if you prefer $W$, lying in an infinite dimensional linear space (the space of density matrices, or the space of Wigner functions; these are just two concrete representations in which the experimenter has particular interest). The data has a probability distribution which is a linear transform of the parameter. Considered as an analytical problem, we have an ill-posed inverse problem, but one which has a lot of beautiful mathematical structure and about which a lot is known (for instance, close connection to the Radon transform). Moreover it has features in common with nonparametric missing data problems (the projections from bivariate to univariate, for instance, and there are more connections we will mention later) and with nonparametric density 
and regression estimation. Thus we think that the time is ripe for this problem to be 
``cracked'' by mathematical and computational statisticians. In this paper we will present some first steps in that direction.

Our main theoretical results are consistency theorems for two estimators. Both estimators are based on approximating the infinite dimensional parameter $\rho$ by a finite dimensional parameter, in fact, thinking of $\rho$ as an infinite dimensional matrix, we simply truncate it to an $N\times N$ matrix where the truncation level $N$ will be allowed to grow with the number of observations $n$. The first estimator employs some analytical inverse formulas expressing the elements of $\rho$ as mean values of certain functions, called pattern functions, of the observations $(X,\Phi)$. Simply replace the theoretical means by empirical averages and one has unbiased estimators of the elements of $\rho$, with moreover finite variance. 
If one applies this technique without truncation the estimate of the matrix $\rho$ as a
 whole will typically not satisfy the nonnegativity constraints. The resulting estimator 
will not be consistent either, with respect to natural distance measures. But provided the truncation level grows with $n$ slowly enough, the overall estimator will be consistent. 
The effect of truncating the density matrix $\rho$ is to project on the subspace generated by the first elements of the corresponding basis, we shall call it the Pattern Function Projection estimator (PFP).

The second estimator we study is sieve maximum likelihood (SML) based on the same truncation to a finite dimensional problem. The truncation level $N$ has to depend on sample size $n$ in order to balance bias and variance. We prove consistency of the SML estimator under an appropriate choice of $N(n)$ by applying a general theorem of \cite{Geer}. To verify the conditions we need to bound certain metric entropy integrals (with bracketing) which express the size of the statistical model under study.

This turns out to be feasible, and indeed to have an elegant solution, by exploiting features of the mapping from parameters (density matrices) to distributions of the data. Various distances between probability distributions possess analogues as distances between density matrices, the mapping from parameter to data turns out to be a contraction, so we can bound metric entropies for the statistical model for the data with quantum metric entropies for the class of density matrices. And the latter can be calculated quite conveniently.

Our results form just a first attempt at studying the statistical properties of estimators 
which are already being used by experimental physicists, but they show that the basic problem is both rich in interesting features and tractable to analysis. The main results so far are consistency theorems for PFP and SML estimators, of both the density matrix and the Wigner function. These results depend on an assumption of the rate at which a truncated density matrix approximates the true one. It seems that the assumption is satisfied for the kinds of states which are met with in practice. However, further work is needed here to describe in physically interpretable terms, when the estimators work. 
Secondly, we need to obtain rates of consistency and to further optimize the construction of the estimator. 
Thirdly, one should explore the properties of penalized maximum likelihood. This will make the truncation level data driven. Fourthly, one should try to make the estimators adaptive with respect to the smoothness of the parameter.
We largely restrict attention to an ideal case of the problem where there is no further noise in  the measurements. In practice, the observations have added to them Gaussian disturbances of known variance. There are some indications that when the variance is larger than a threshold 
of $1/2$, reconstruction becomes qualitatively much more difficult. 
This needs to be researched from the point of view of optimal rates of convergence.

We also only considered one particular though quite convenient way of sieving the model, i.e., one particular class of finite dimensional approximations. There are many other possibilities and some of them might allow easier analysis and easier computation. For instance, instead of truncating the matrix $\rho$ in a given basis, one could truncate in an arbitrary basis, so that the finite dimensional approximations would correspond to specifying $N$ arbitrary state vectors and a probability distribution over these ``pure states''. Now the problem has become a missing data problem, where the ``full data'' would assign to each observation also the label of the pure state from which it came. In the full data problem we need to reconstruct not a matrix but a set of vectors, together with an ordinary probability distribution over the set, so the ``full data'' problem is statistically speaking a much easier problem that the missing data problem. We shall use a version of this to derive an Expectation-Maximization algorithm (EM) for the practical implementation of the SML estimator, see Section \ref{sec.experimental}. One could imagine that Bayesian reconstruction methods could also exploit this structure.
 
The analogy with density estimation could suggest new statistical approaches here. 
Finally, it is most important to add to the estimated parameter, estimates of its accuracy. This is absolutely vital for applications, but so far no valid approach is available.

In Section \ref{sec.problem} we introduce first, very briefly, the statistical problems we are concerned with. We then give a short review of the basic notions of quantum mechanics which are needed in this paper. Concepts such as observables, states, measurements and quantum state tomography are explained by using finite dimensional complex matrices. At the end of the section we show how to generalize this to the infinite dimensional case and describe the experimental set-up of Quantum Homodyne Detection pointing out the relation with computerized tomography.

In Section \ref{sec.DensityMatrixEstimation} we present results on consistency of density matrix estimators: projection estimator based on pattern functions, and sieve maximum likelihood estimator. The last subsection extends the results of the previous ones to estimating the Wigner function. 

Section \ref{sec.efficiency} deals with the detection losses occurring in experiments due in part to the inefficiency of the detectors. This adds a deconvolution problem on the top of our tomography estimation.

Section \ref{sec.experimental} shows experimental results. We illustrate the behavior of the studied estimators and propose some practical tools for the implementation --- EM algorithm. The last section finishes with some concluding remarks to the whole paper and open problems. The main purpose of the paper is to bring the attention of the statistical community to these problem; thus, some proofs are just sketched. They fully appear in a complementary paper, \cite{Guta}.

\nopagebreak

\section{Physical background}
\label{sec.problem}

Our statistical problem is to reconstruct the quantum state of light by using 
data obtained from measuring identical pulses of light through a technique called 
Quantum Homodyne Detection (QHD). In particular we will estimate the quantum state in two different representations or parameterizations: the density matrix and the Wigner function. The physics behind this statistical problem is presented in subsection \ref{sec.qmeasurements} which serves as introduction to basic notions of quantum theory, and subsection \ref{sec.qhomodyne} which describes the model of quantum homodyne detection from quantum optics. The relations between the different parameters of the problem are summarized in the diagram at the end of subsection \ref{sec.qhomodyne} followed by Table \ref{tbl.states} containing some examples of quantum states. 
For the reader who is not interested in the quantum physics background we state the statistical problem in the next subsection and we will return to it in Section \ref{sec.DensityMatrixEstimation}. 

\subsection{Summary of statistical problem}

We observe $(X_1,\Phi_1),\dots ,(X_n,\Phi_n)$, i.i.d.~random variables with values in 
$\mathbb{R}\times [0,\pi]$ and  distribution $P_\rho$ depending on 
the unknown parameter $\rho$ which is an infinite dimensional matrix 
$\rho = [\rho_{j,k}]_{j,k=0,\ldots,\infty}$ such that $\mathrm{Tr} \rho =1$ 
(trace one) and $\rho \geq 0$ (positive definite). The probability density of $P_\rho$ is
\begin{equation}\label{eq.concrete.p}
p_\rho(x,\phi)=\frac{1}{\pi}\sum_{j,k=0}^\infty \rho_{j,k}\psi_k(x)\psi_j(x)e^{-i(j-k)\phi},
\end{equation}
where the functions $\{\psi_n\}$ to be specified later, form a orthonormal basis of the 
space of complex square integrable functions on $\mathbb{R}$. Because $\rho$ is positive definite and has trace $1$, this is a probability density: real, nonnegative, integrates to $1$.  
The data $(X_\ell,\Phi_\ell)$ comes from independent QHD measurements on identically prepared pulses of light whose properties or state are completely encoded in the matrix $\rho$ called a density matrix. We will consider the problem of estimating $\rho$ from a given sample. 

Previous attempts by physicists to estimate the density matrix $\rho$ have focused mainly on the estimation of the individual matrix elements without considering the accuracy of the estimated density matrix with respect to natural distances of the underlying parameter space. In Section \ref{sec.DensityMatrixEstimation} we will present 
consistency results in the space of density matrices w.r.t. $L_1$ and $L_2$-norms 
using two different types of estimators, namely projection and sieve maximum likelihood estimators.

An alternative representation of the quantum state is through the Wigner function 
$W_{\rho}:\mathbb{R}^2 \to\mathbb{R}$ whose estimation is close to a classical 
computer tomography problem namely, Positron Emission Tomography (PET), 
\cite{Vardi}. In PET one would like to estimate a probability density $f$ on $\mathbb{R}^2$ from i.i.d.~observations $(X_1,\Phi_1), \dots, (X_n,\Phi_n)$, with probability density equal to the Radon transform of $f$:  
\begin{equation*}
\mathcal{R}[f](x,\phi)=\int_{-\infty}^\infty 
f(x\cos\phi+t\sin\phi, x\sin\phi-t\cos\phi)dt.
\end{equation*}
Although the Wigner function is in general {\it not} positive, its Radon transform is always a probability density, in fact $\mathcal{R}[W_{\rho}](x,\phi) = p_\rho(x,\phi)$. 
As the Wigner function $W_\rho$ is in one-to-one correspondence with the density matrix $\rho$, our state reconstruction problem can be stated as to estimate the Wigner function $W_{\rho}$. This is an ill posed inverse problem as seen from the formula for the inverse of the Radon transform
\begin{equation}\label{eq.inverse.Radon.transform}
W_{\rho}(q,p) = \frac{1}{2 \pi^2} \int_{0}^{\pi} \int_{-\infty}^{+\infty} p_\rho(x,\phi) 
K(q\cos\phi+p\sin\phi-x) \, dx d\phi,
\end{equation}
where
\begin{equation}\label{eq.kernel}
K(x) = \frac{1}{2} \int_{-\infty}^{+\infty} |\xi| \exp(i\xi x) d\xi,
\end{equation}
makes sense only as a generalized function. To correct this one usually makes a cut-off in the range of the above integral and gets a well behaved kernel function 
$K_c(x)=\frac{1}{2} \int_{-c}^{c} |\xi| \exp(i\xi x) d\xi$. 
Then the tomographic estimator of $W_{\rho}$ is the average sampled kernel
\begin{equation*}
\widehat{W}_{c,n}(q,p) = \frac{1}{2 \pi^2n} \sum_{\ell=1}^{n} 
K_c(q\cos\Phi_{\ell}+p\sin\Phi_{\ell}-X_\ell).
\end{equation*}
For consistency one needs to let the `bandwidth' $h=1/c$ depend on the sample size $n$ and $h_n\to 0$ as $n\to\infty$ at an appropriate rate. 

In this paper we will not follow this approach, which will be treated separately in future work. Instead, we use a plug-in type estimator based on the property
\begin{equation} \label{eq.linearity.Wigner.rho}
W_{\rho}(q,p) = \sum_{k,j=0}^{\infty} \rho_{k,j} W_{k,j}(q,p),
\end{equation}
where $W_{k,j}$'s are known functions and we replace $\rho$ by its above mentioned estimators. We shall prove consistency of the proposed estimators of the Wigner function w.\ r.\ t.~$L_2$ and supremum norms in the corresponding space.

\subsection{Quantum systems and measurements}
\label{sec.qmeasurements}

This subsection serves as a short introduction to the basic notions of quantum mechanics which will be needed in this paper. For simplicity we will deal first with finite 
dimensional quantum systems and leave the infinite dimensional case for the next subsection. For further details on quantum statistical inference 
we refer to the review \cite{Barndorff-Nielsen&Gill&Jupp} and the classic 
textbooks \cite{Helstrom} and \cite{Holevo}. 

In classical mechanics the state of macroscopic systems like billiard balls, pendulums or stellar systems is described by points on a manifold or ``phase space'', each of the point's coordinates corresponding to an attribute which we can measure such as position and momentum. 
Therefore the functions on the phase space are called observables. When there exists 
uncertainty about the exact point in the phase space, or we deal with a statistical 
ensemble, the state is modelled by a probability distribution on the phase space, 
and the observables become random variables. 

Quantum mechanics also deals with observables such as position and momentum of a particle, spin of an electron, number of photons in a cavity, but breaks from classical mechanics in that these are no longer represented by functions on a phase space but by Hermitian matrices, that is, complex valued matrices which are invariant under transposition followed by complex conjugation. For example, the components in different directions of the spin of an electron are certain $2\times 2$ complex Hermitian 
matrices $\sigma_x,\sigma_y, \sigma_z$.

Any $d$-dimensional complex Hermitian matrix $\mathbf{X}$ can be diagonalized by changing the standard basis of $\mathbb{C}^d$ to another orthonormal basis $\{e_1,\dots, e_d\}$ such that $\mathbf{X}e_i=x_i e_i$ for $i=1,\dots, d$, with $x_i\in \mathbb{R}$. The vectors $e_i$ and numbers $x_i$ are called 
eigenvectors and respectively eigenvalues of $\mathbf{X}$. With respect to the new basis we can write    
\begin{displaymath}\label{eq.diagonalization}
\mathbf{X}= 
\left( \begin{array}{ccccc}
 x_1 &  0  & 0   & \ldots &  0 \\ 
  0  & x_2 & 0   & \ldots &  0 \\
  0  &  0  & x_3 & \ldots &  0 \\
     &     &     & \ddots &    \\
  0  &  0  &  0  & \ldots & x_d
\end{array}\right)
\end{displaymath}
The physical interpretation of the eigenvalues is that when measuring the 
observable $\mathbf{X}$ we obtain (randomly) one of the values $x_i$ according to a probability distribution depending on the state of the system before measurement and on the observable $\mathbf{X}$. This probability measure is degenerate if and only if the system before measurement was prepared in a special state called an eigenstate of $\mathbf{X}$. We represent such a state mathematically by the projection $\mathbf{P}_i$ onto the one dimensional space generated by the vector $e_i$ in $\mathbb{C}^d$. Given a probability distribution $\{p_1,\dots, p_d\}$ over the finite set $\{x_1,\dots, x_d\}$, we describe a statistical ensemble in which a proportion $p_i$ of systems is prepared in the state $\mathbf{P}_i$ by the convex combination $\rho=\sum_i p_i \mathbf{P}_i$. The expected value of the random result $X$ when measuring the observable $\mathbf{X}$ for this particular state is equal to $\sum_i p_i x_i$ which can be written shortly
\begin{equation}\label{eq.quantumexpectation}
\mathbb{E}_\rho(X):=\mathrm{Tr}(\rho \mathbf{X}).
\end{equation}
Similarly, the probability distribution can be recovered as 
\begin{equation}\label{eq.quantumprobability}
p_i=\mathrm{Tr}(\rho \mathbf{P}_i)
\end{equation}
thanks to the orthogonality property 
$\mathrm{Tr}( \mathbf{P}_i\mathbf{P}_j)=\delta_{ij}$. 

Now, let $\mathbf{Y}$ be a different observable and suppose that $\mathbf{Y}$ does not 
commute with $\mathbf{X}$, that is $\mathbf{X}\mathbf{Y}\neq \mathbf{Y}\mathbf{X}$, then the two observables cannot be diagonalized in the same basis, their eigenvectors are different. Consequently, states which are mixtures of eigenvectors of $\mathbf{X}$ typically will not be mixtures of eigenvectors of $\mathbf{Y}$ and vice-versa. This leads to an expanded formulation of the notion of state in 
quantum mechanics independent of any basis associated to a particular observable, 
and the recipe for calculating expectations and distributions of measurement 
results. 

Any preparation procedure results in an statistical 
ensemble, or state, which is described mathematically by a matrix $\rho$ with the 
following properties 
\begin{enumerate}
\item
 $\rho \geq 0$  (positive definite matrix),
\item
 $\mathrm{Tr}(\rho)=1$ (normalization).
\end{enumerate}
In physics $ \rho$ is called a {\it density matrix}, and is for a quantum mechanical system an analogue of a probability density. Notice that the special state 
$\sum_i p_i \mathbf{P}_i$ defined above is a particular case of density matrix, since it is a mixture of eigenstates of the observable $\mathbf{X}$. The density matrices of dimension $d$ form a convex set $\mathcal{S}_d$, whose extremals are the {\it pure} or {\it vector} states, represented by orthogonal projections $\mathbf{P}(\psi)$ onto one 
dimensional spaces spanned by {\it arbitrary} vectors $\psi\in\mathbb{C}^d$. Any state can be represented as a mixture of pure states which are not necessarily eigenstates of a particular observable.

When measuring an observable, for example $\mathbf{X}$, of a quantum system prepared in the state $\rho$ we obtain a random result $X\in\{x_1,\dots ,x_d\}$ with 
probability distribution given by equation (\ref{eq.quantumprobability}), 
expectation as in equation (\ref{eq.quantumexpectation}), and characteristic function 
\begin{equation}\label{eq.generating.function}
G(t):=\mathbb{E}_\rho\big(\exp(itX)\big)=\mathrm{Tr}\big(\rho\exp(it\mathbf{X})\big).
\end{equation}
In order to avoid confusion we stress the important 
difference between $\mathbf{X}$ which is a matrix and $X$ which is a real-valued random variable. More concretely, if we write $\rho$ in the basis of eigenvectors of $\mathbf{X}$ then we 
obtain the map $\mathbf{M}$ from states to probability distributions $ P_\rho^{(\mathbf{X})}$ 
over results $\{x_1,\dots , x_d\}$
\begin{displaymath}\label{eq.measurement.map}
\mathbf{M}:\rho= 
\left( \begin{array}{cccc}
 \rho_{11} &  \rho_{12}   &  \ldots &  \rho_{1d}   \\ 
 \rho_{21} &  \rho_{22}   &  \ldots &  \rho_{2d}   \\
           &              &  \ddots &              \\
 \rho_{d1} &  \rho_{d2}   &  \ldots & \rho_{dd} 
\end{array}\right)
\longmapsto 
P_\rho^{(\mathbf{X})}=
\left( \begin{array}{cccc}
 \rho_{11}\\
 \rho_{22}\\
 \vdots   \\
 \rho_{dd}\\
\end{array}\right).
\end{displaymath}
Notice that $P_\rho^{(\mathbf{X})}$ is indeed a probability distribution as a consequence of the defining properties of states, and it does not contain information about the off-diagonal elements of $\rho$, meaning that measuring only the observable $\mathbf{X}$ is not enough to identify the unknown state. Roughly speaking, as 
$\mathrm{dim}(\mathcal{S}_d)=d^2-1=(d-1)(d+1)$, one has to measure on many 
identical systems each one of a number of $d+1$ mutually non-commuting observables in order to have a one-to-one map between states and probability distributions of results. The probing of identically prepared quantum systems from different `angles' in order to reconstruct their state is broadly named {\it quantum state tomography} in the physics literature. 

Let us suppose that we have at our disposal $n$ systems identically prepared in an unknown state $\rho\in\mathcal{S}_d$, and for each of the systems we can measure one of the fixed observables $\mathbf{X}(1),\dots ,\mathbf{X}(d+1)$. We write the observables in diagonal form
\begin{equation}
\mathbf{X}(i)=\sum_{a=1}^d x_{i,a}\mathbf{P}_{i,a}
\end{equation}
where $x_{i,a}$ eigenvalues and $\mathbf{P}_{i,a}$ eigenstates. We will perform a randomized experiment, i.e.~for each system we will choose the observable to be measured by randomly selecting its index according to a probability distribution 
$P^\mathbf{\Phi}$ over $\{1, \dots ,d+1\}$. 
The results of the measurement on the $k^{\text{th}}$ system are the pair $Y_k=(X_k, \Phi_k)$
where $\Phi_1,\dots ,\Phi_n$ are i.i.d. with probability distribution $P^{(\Phi)}$ and 
$X_k$ is the result of measuring the observable $\mathbf{X}(\Phi_k)$ whose conditional 
distribution is given by 
\begin{equation}
P^\mathbf{M}(X_k=x_{i,a}|\Phi_k=i)=\mathrm{Tr}(\rho \mathbf{P}_{i,a})
\end{equation}
The statistical problem is now to estimate the parameter $\rho$ from the data 
$Y_1,Y_2,\dots ,Y_n$. In the next subsection we will describe quantum homodyne tomography as an analogue of this problem for infinite dimensional systems.

\subsection{Quantum homodyne tomography}
\label{sec.qhomodyne}

Although correct and sufficient when describing certain quantum properties such as the 
spin of a particle, the model presented above needs to be enlarged in order to cope with 
quantum systems with `continuous variables' which will be central in our statistical problem. This technical point can be summarized as follows: we replace $\mathbb{C}^d$ by an infinite dimensional complex Hilbert space $\mathcal{H}$, the Hermitian matrices becoming {\em selfadjoint operators} acting on $\mathcal{H}$. The spectral theorem tells us that selfadjoint operators can be `diagonalized' in the spirit of 
(\ref{eq.diagonalization}) but the spectrum (the set of `eigenvalues') can have a more complicated structure, for example it can be continuous as we will see below. 
The density matrices are {\em positive} selfadjoint operators $\mathbf{\rho}$ such that $\mathrm{Tr}(\rho)=1$ and can be regarded as infinite dimensional matrices with elements $\rho_{k,j} := \br \psi_j, \rho \psi_k \ke$ for a given orthonormal basis $\{\psi_1,\psi_2,\dots\}$ in $\mathcal{H}$. 
 
The central example of a system with continuous variables in this paper is the quantum particle. Its basic observables position and momentum, are two unbounded selfadjoint operators $\mathbf{Q}$ and $\mathbf{P}$ respectively, acting on $L^2(\R)$, the space of square integrable complex valued functions on $\R$
\begin{eqnarray} 
&&(\mathbf{Q}\psi_1)(x) = x\psi_1(x),\nonumber\\
&&(\mathbf{P}\psi_2)(x) = -i\frac{d\psi_2(x)}{dx},\nonumber
\end{eqnarray}
for $\psi_1,\psi_2$ arbitrary functions. The operators satisfy 
Heisenberg's {\it commutation relations} 
$\mathbf{Q}\mathbf{P}-\mathbf{P}\mathbf{Q}=i\mathbf{1}$ which implies that they cannot be measured simultaneously. The problem of (separately) measuring such observables has been elusive until ten years ago when pioneering experiments in quantum optics by \cite{Smithey}, led to a powerful measurement technique called {\it quantum homodyne detection}. This technique is the basis of a continuous analogue of the measurement scheme presented at the end of the previous subsection where $d+1$ observables were measured in the case of a $d$-dimensional quantum system.

\begin{figure}[ht]
\begin{center}
\includegraphics[width=8cm]{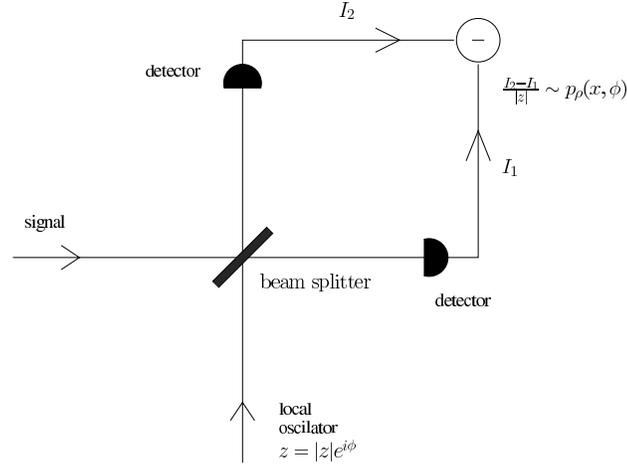}
\caption{Quantum Homodyne Tomography measurement system}\label{fig.QHtomography}
\end{center}
\end{figure}

The quantum system to be measured is a beam of light with a fixed frequency whose observables are the electric and magnetic field amplitudes which satisfy commutation relations identical to those characterizing the quantum particle, with which they will be identified from now on. Their linear combinations 
$\mathbf{X}_\phi=\cos \phi\mathbf{Q}+\sin\phi \mathbf{P}$ are called {\it quadratures}, 
and homodyne detection is about measuring the quadratures for {\it all} phases 
$\phi\in [0,\pi]$. The experimental setup shown in Figure \ref{fig.QHtomography} 
contains an additional laser called a local oscillator (LO) of high intensity $|z|\gg 1$ and relative phase $\phi$ with respect to the beam in the unknown state $\rho$. The two beams are combined through a fifty-fifty beam splitter, and the two emerging beams are then measured by two photon detectors. A simple quantum optics computation (see \cite{Leonhardt}) shows that in the limit of big LO intensity the difference of the measurement results (countings) of the two detectors rescaled by the LO intensity $X=\frac{I_1-I_2}{|z|}$ has the probability distribution corresponding to the measurement of the quadrature $\mathbf{X}_\phi$. 
The result $X$ takes values in $\R$ and its probability distribution $P_\rho(\cdot|\phi)$ 
has a density $p_\rho(x|\phi)$ and characteristic function
(see equation (\ref{eq.generating.function})) 
\begin{equation}\label{eq.generating_function}
G_{\rho,\mathbf{X}_\phi}(t)=\mathrm{Tr}\big(\rho \exp(it\mathbf{X}_\phi)\big). 
\end{equation}
The phase $\phi$ can be controlled by the experimenter by adjusting a parameter of the local oscillator. We assume that he chooses it randomly uniformly distributed over the 
interval $[0,\pi]$. Then the joint probability distribution $P_\rho$ for the pair consisting in 
measurement result and phase $Y:=(X,\Phi)$ has density $p_\rho(x,\phi)$ equal to 
$\frac{1}{\pi}p_\rho(x|\phi)$ with respect to the measure $dx\times d\phi$ on $\R\times[0,\pi]$. 
An attractive feature of the homodyne detection scheme is the invertibility of the map 
$\mathbf{T}$ that associates $P_\rho$ to $\rho$, making it possible to asymptotically infer the unknown 
parameter $\rho$ from the i.i.d. results $Y_1,\dots, Y_n$ of homodyne measurements on $n$ systems prepared in the state $\rho$.

We will see now why this state estimation method is called {\it quantum homodyne tomography} by drawing a parallel with computerized tomography used in the hospitals. 
In quantum optics it is common to represent the state of a quantum system by a 
certain function on $\R^2$ called the {\it Wigner function} $W_{\rho}(q,p)$ which is much like a joint probability density for $\mathbf{Q}$ and $\mathbf{P}$, for instance its 
marginals are the probability densities for measuring $\mathbf{Q}$ and respectively 
$\mathbf{P}$. The Wigner function of the state $\rho$ is defined by demanding that its Fourier transform $\mathcal{F}_2$ with respect to {\it both} variables has the following property
\begin{equation}\label{def.Wigner}
\widetilde{W}_\rho(u,v):=\mathcal{F}_2[W_{\rho}](u,v)=
\mathrm{Tr}\big(\rho \exp(-iu\mathbf{Q}-iv\mathbf{P})\big).
\end{equation}
We see from this equation that if $\mathbf{Q}$ and $\mathbf{P}$ were commuting operators then $W(q,p)$ would indeed be the joint probability density of outcomes of their measurement. As the two observables cannot be measured simultaneously, we cannot speak of a joint density, in fact the Wigner function need not be positive, but many interesting features of the quantum state can be visualized in this way. Let $(u,v)=(t\cos\phi, t\sin\phi)$, then 
\begin{equation}\label{eq.Fourier transform.tildeW}
\widetilde{W}(u,v)=\mathcal{F}_1[p(\cdot,\phi)](t)=
\mathrm{Tr}\big(\rho \exp(-it\mathbf{X}_\phi)\big)
\end{equation}
where the Fourier transform $\mathcal{F}_1$ in the last term is with respect to the first variable, keeping $\phi$ fixed. 

The equations (\ref{def.Wigner}) and (\ref{eq.Fourier transform.tildeW}) are well known in the theory of Radon transform $\mathcal{R}$ and imply that for each fixed $\phi$, the 
probability density $p_\rho(x,\phi)$ is the marginal of the Wigner function with respect to the direction $\phi$ in the plane, 
\begin{equation}
p_\rho(x,\phi)=\mathcal{R}[W_\rho](x,\phi)=\int_{-\infty}^\infty 
W_{\rho}(x\cos\phi+y\sin\phi, x\sin\phi-y\cos\phi)dy,
\end{equation}
adding quantum homodyne tomography to a number of applications ranging from computerized tomography to astronomy and geophysics, \cite{Deans}. 
In computerized tomography one reconstructs an image of the tissue distribution in a cross-section of the human body by recording events whereby pairs of positrons emitted by an injected radioactive substance hit detectors placed in a ring around the body after 
flying in opposite directions along an axis determined by an angle $\phi\in [0,\pi]$. In quantum homodyne tomography the role of the unknown distribution is played by 
the Wigner function which is in general not positive, but has a probability density
$p_\rho(x|\phi)$ as marginal along any direction $\phi$.

The following diagram summarizes the relations between the various objects in our problem: \\

\hspace{0.1cm}
\begindc[3]
\obj(10,30){$\rho$}
\obj(30,30){$W_{\rho}$}
\obj(65,30){$p_{\rho}(x,\phi)$}
\obj(110,30){$(X_{1},\Phi_{1}),\dots, (X_n,\Phi_n)$.}
\obj(47,10){$\widetilde{W}_{\rho}$}
\mor(11,30)(29,30){}[\atleft,\solidarrow]
\mor(29,30)(11,30){}[\atleft,\solidarrow]
\mor(31,30)(60,30){$\mathcal{R}$}[\atleft,\solidarrow]
\mor(31,29)(45,11){$\mathcal{F}_{2}$}[\atright,\solidarrow]
\mor(60,29)(47,11){$\ \mathcal{F}_{1}$}[\atleft,\solidarrow]
\mor(70,30)(95,30){experiment}[\atleft,\solidarrow]
\enddc

\vspace{2mm}

Finally in Table \ref{tbl.states} we give some examples of density matrices and their 
corresponding Wigner function representations for different states. The matrix elements 
$\rho_{k,j}$ are calculated with respect to the orthonormal base corresponding to the wave functions of $k$ photons states
\begin{equation}\label{eq.psi_n}
\psi_k(x)=H_k(x)e^{-x^2/2}
\end{equation} 
where $H_k$ are the Hermite polynomials normalized such that $\int \psi_k^2=1$.  
A few graphical representations can be seen in Figure \ref{fig.examples}. 

\begin{longtable}{|c|c|c|}
\hline
State & Density matrix $\rho_{k,j}$ & Wigner function $W(q,p)$\\
\hline 
& & \\[-1mm]
Vacuum state 
& $\rho_{0,0}=1,\qquad\mathrm{rest~zero}$
& $\frac{1}{\pi} \exp (-q^2-p^2)$ 
\\[2mm]
Single Photon state 
& $\rho_{1,1}=1, \qquad\mathrm{rest~zero}$
& $\frac{1}{\pi} (2q^2+2p^2-1) \exp (-q^2-p^2)$ 
\\[2mm]
Thermal state $\beta > 0$ 
& $ \delta_k^j (1-e^{-\beta}) e^{-\beta k}$ & $\frac{1}{\pi} 
\tanh(\beta/2) \exp [-(q^2+p^2)\tanh(\beta/2)]$ \\[2mm]
Coherent state, $N \in \mathbb{R}_+$ 
& $\exp(-N) \frac{N^{k+j}}{\sqrt{k!j!}}$ & $\frac{1}{\pi} 
\exp (-(q-\sqrt{N})^2-p^2)$ \\[2mm]
Squeezed  state
&$C(N,\xi)(\frac{1}{2}\tanh(\xi))^{k+j}  \times  $ 
& $\frac{1}{\pi} \exp (-e^{2 \xi}(q-\alpha)^2- e^{-2 \xi} p^2)$ \\[1mm]
$N \in \mathbb{R}_+$ ,$\quad \xi \in \mathbb{R},$ 
& $ H_{j}(\gamma)H_{k}(\gamma)/\sqrt{j!k!} $ 
& \\[1mm]
\hline
\caption{Density matrix and Wigner function of some quantum states}\label{tbl.states}
\end{longtable}

The vacuum is the pure state of zero photons, notice that in this case the distributions of $\mathbf{Q}$ and $\mathbf{P}$ are Gaussian. The thermal state is a mixed 
state describing equilibrium at temperature $T=1/\beta$, having Gaussian Wigner function with variance increasing with the temperature. The coherent state is pure and 
characterizes the laser pulse. The photon number is Poisson distributed with an average of $N$ photons. The squeezed states have Gaussian Wigner functions whose variances for the two directions are different but have a fixed product.  
The parameters $N$ and $\xi$ satisfy the condition $N\geq\sinh^2(\xi)$, $C(N,\xi)$ is a normalization constant, $\alpha=\frac{(N-\sinh^2(\xi))^{1/2}}{\cosh(\xi)-\sinh(\xi)}$, and 
$\gamma = (\frac{\alpha}{\sinh(2\xi)})^{1/2}$.  

\section{Density matrix estimation}
\label{sec.DensityMatrixEstimation}

\cite{D'Ariano.0} presented the density matrix analogue of formula 
(\ref{eq.inverse.Radon.transform}) of the Wigner function as inverse Radon transform  of the probability density $p_\rho$
\begin{equation}\label{eq.inverseT}
\rho=\int_{-\infty}^{\infty}dx \int_0^\pi \frac{d\phi}{\pi}
p_\rho(x,\phi)K(\mathbf{X}_\phi-x \mathbf{1}),
\end{equation}
where $K$ is the generalized function given in equation (\ref{eq.kernel}) whose argument is a selfadjoint operator $\mathbf{X}_\phi-x\mathbf{1}$. 
The method has been further analyzed in \cite{D'Ariano.1}, \cite{D'Ariano.3}, see also \cite{D'Ariano.2}. We recall that in the case of the Wigner function we needed to regularize the kernel $K$ by introducing a cut-off in the integral (\ref{eq.kernel}). For density matrices the philosophy will be rather to project on a finite dimensional subspace of $L_2(\R)$ whose dimension $N$ will play the role of the cut-off. In fact all the matrix elements of the density matrix $\rho$ with respect to the orthonormal basis $\{\psi_k\}_{k=0}^\infty$ defined in (\ref{eq.psi_n}), can be expressed as kernel integrals
\begin{equation}\label{eq.quantum_tomographic_rho_n,n+d}
\rho_{k,j}=\int_{-\infty}^\infty dx \int_0^\pi \frac{d\phi}{\pi} p_\rho(x,\phi)
f_{k,j}(x)e^{-i(j-k)\phi}, 
\end{equation}
with $f_{k,j}=f_{j,k}$ bounded real valued functions which in the quantum tomography literature are called {\it pattern functions}. The singularity of the kernel $K$ is reflected in the asymptotic behavior of $f_{k,j}$ as $k,j\to\infty$. 
A first formula for $f_{k,j}$ was found in \cite{D'Ariano.3} and uses Laguerre polynomials. This was followed by a more transparent one due to \cite{Leonhardt.Richter},
\begin{equation}\label{eq.pattern_functions_leonhardt}
f_{k,j}(x)=\frac{d}{dx}(\psi_k(x)\varphi_{j}(x)),
\end{equation}
for $j \geq k$, where $\psi_k$ and $\varphi_j$ represent the square integrable and respectively the unbounded solutions of the Schr\"{o}dinger equation, \begin{equation}\label{eq.Schrodinger}
\left[-\frac{1}{2}\frac{d^2}{dx^2}+\frac{1}{2}x^2\right]\psi=\omega~\psi ,
\qquad \omega\in\R.
\end{equation}
\begin{figure}[h]
\begin{center}
\includegraphics[width=10cm]{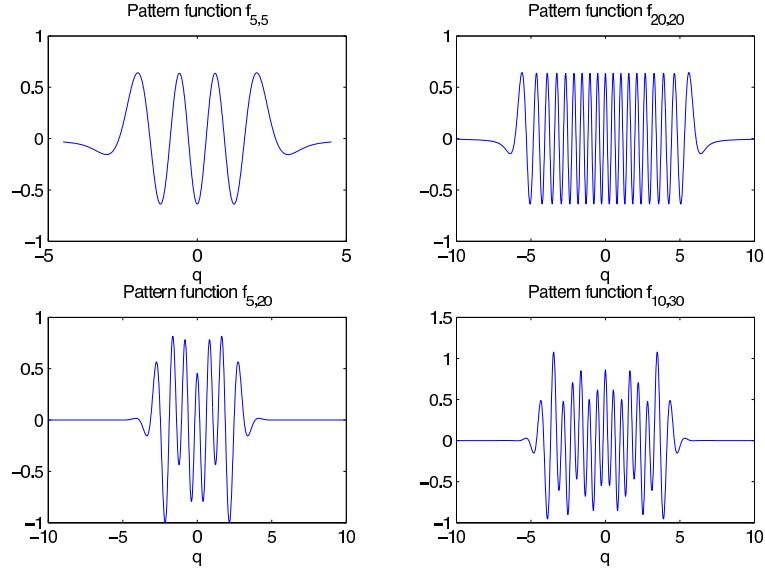}
\caption{Pattern functions}\label{fig.ptnfnc}

\end{center}
\end{figure}
Figure \ref{fig.ptnfnc} shows pattern functions for different values of $k$ and $j$. 
We notice that the oscillatory part is concentrated in an interval centered at zero 
whose length increase with $k$ and $j$, the number of oscillations increases with 
$k$ and $j$ and the functions become more irregular as we move away from the diagonal. It can be shown that tails of the pattern function decay like $x^{-2- |k-j|}$. 
More properties of the pattern function can be found in \cite{Leonhardt.Richter} and \cite{Guta}.

\subsection{Pattern function projection estimation}
\label{sec.pfe}

Equation (\ref{eq.quantum_tomographic_rho_n,n+d}) suggests the {\it unbiased estimator} $\hat{\rho}^{(n)}$ of $\rho$, based on $n$ i.i.d. observations of $(X,\Phi)$, 
whose matrix elements are:
\begin{equation}\label{eq.unbiased_estimator}
\hat{\rho}^{(n)}_{k,j}=\frac{1}{n}\sum_{\ell=1}^n F_{k,j}(X_\ell,\Phi_\ell),
\end{equation} 
where $F_{k,j}(x,\phi) = f_{k,j}(x) e^{-i(j-k)\phi}$, see 
\cite{D'Ariano.2}, \cite{D'Ariano.3, Leonhardt.Richter}. By the strong law of large numbers the 
individual matrix elements of this estimator converge to the matrix elements of the true
 parameter $\rho$. However the infinite matrix $\hat{\rho}^{(n)}$ need not be positive, 
normalized, or even selfadjoint, thus it cannot be interpreted as a quantum state. 
These problems are similar to those encountered when trying to estimate an unknown probability density by using unbiased estimators for all its Fourier coefficients. The remedy is to estimate only a finite number of coefficients at any moment, obtaining a projection estimator onto the subspace generated by linear combinations of a finite subset of the basis vectors. In our case we will project onto the space of matrices of dimension $N=N(n)$ with respect to the basis $\{\psi_k\}_{k=0}^\infty$,
\begin{displaymath}
\hat{\rho}^{(N,n)}_{k,j}=\frac{1}{n}\sum_{\ell=1}^n F_{k,j}(X_\ell,\Phi_\ell),
\qquad \mathrm{for}~ 0\leq k,j\leq N-1,
\end{displaymath}
and $\hat{\rho}^{(n)}_{k,j}=0$ for $ (k\wedge j) \geq N$.

In order to test the performance of our estimators we introduce the $L_1$ and $L_2$ distances 
on the space of density matrices. Let $\rho$ and $\tau$ be two density matrices with $\rho-\tau=\sum_i \lambda_i\mathbf{P}_i$ the diagonal form of their difference, and notice that some of the eigenvalues are positive and some negative such that their sum is zero due to the normalization of the density matrices. 
We define the absolute value $|\rho-\tau|:=\sum_i |\lambda_i|\mathbf{P}_i$ and the norms  
\begin{eqnarray}
&&\|\rho-\tau\|_1:=\mathrm{Tr}(|\rho-\tau|),\\
&&\|\rho-\tau\|_2:=\left[\mathrm{Tr}(|\rho-\tau|^2)\right]^{1/2}=
                 \left[\sum_{j,k= 0}^\infty|\rho_{k,j}-\tau_{k,j}|^2\right]^{1/2}.
\end{eqnarray}

Let us consider now the mean integrated square error (MISE) and split it into the bias and variance parts:
\begin{equation}
\mathbb{E}\big(\sum_{j,k=0}^\infty|\hat{\rho}^{(N,n)}_{k,j}-\rho_{k,j}|^2\big)=
\sum_{(k\wedge j)\geq N}|\rho_{k,j}|^2  + 
\mathbb{E}\big(\sum_{j,k= 0}^{N-1} |\hat{\rho}^{(N,n)}_{k,j}-\rho_{k,j}|^2\big)=b^2(n)+\sigma^2(n).
\end{equation} 
By choosing $N(n)\to\infty$ as $n\to\infty$ the bias $b^2(n)$ converges to zero. 
For the variance we have the upper bound
\begin{eqnarray}
\sigma^2(n)&=&
\frac{1}{n}\sum_{j,k= 0}^{N-1} 
\int_0^\pi \int_\R |F_{k,j}(x,\phi)-\rho_{k,j}|^2 p_\rho(x,\phi)d\phi dx \nonumber\\
&\leq& \frac{1}{n} \sum_{j,k= 0}^{N-1}
\int_0^\pi \int_\R |F_{k,j}(x,\phi)|^2 p_\rho(x,\phi)d\phi dx \nonumber\\
&= &\frac{1}{n}\sum_{j,k= 0}^{N-1} \int_\R dx f_{k,j}(x)^2\int_0^\pi p_\rho(x,\phi)d\phi.
\label{eq.variance.upper.bound}
\end{eqnarray} 
The proof of the following lemma on the norms of the pattern functions can be found in \cite{Guta}. 
\begin{lemma} There exist constants $C_1,C_2, C_3$ such that
\begin{eqnarray}\label{eq.sum.pattern.functions}
&& \sum_{k,j\geq0}^N\|f_{k,j}\|_\infty^2 \leq C_1 N^{7/3},\\
&& \sum_{k,j\geq0}^N\|f_{k,j}\|_2^2 \leq C_2 N^3,\\
&& \int_0^{2\pi} p_\rho(x,\phi)d\phi\leq C_3, \qquad \mathrm{for~all}~ x\in\R.
\end{eqnarray}
\end{lemma}
By applying the lemma to equation (\ref{eq.variance.upper.bound}) we conclude that the estimator 
$\hat{\rho}^{(N,n)}$ is consistent with respect to the $L_2$ distance if we choose 
$N(n)\to \infty$ as $n\to \infty$ such that $N(n)=o(n^{3/7})$. Based on the property 
$\|\rho^{(N,n)}\|_1 \leq \sqrt{N} \|\rho^{(N,n)}\|_2$ we can prove a similar result concerning 
$L_1$-consistency, see \cite{Guta}. 

\begin{theorem}\label{th.norm2} Let $N(n)\to \infty$ be the dimension of the pattern function 
projection estimator. If $N(n)=o(n^{3/7})$, then
\begin{equation}
\mathbb{E}\big( \|\hat{\rho}^{(N,n)}-\rho\|_2^2\big)\to 0 ,\qquad \mathrm{as}~ n\to \infty.
\end{equation}
 If $N(n)=o(n^{3/10})$ then 
 \begin{equation}
 \mathbb{E}\big( \|\hat{\rho}^{(N,n)}-\rho\|_1\big)\to 0 ,\qquad \mathrm{as}~ n\to \infty.
 \end{equation}
\end{theorem}
Rates of consistency can be obtained by assuming that the state belongs to a given class for which upper bounds of the bias can be calculated and $N(n)$ is chosen such as to balance bias and variance. This problem will be attacked in future work within the minimax framework. In Section \ref{sec.experimental} we present a data-dependent way of selecting the dimension of the projection estimator based on the minimization of the empirical $L_2$-risk using a cross-validation technique.

\subsection{Sieve maximum likelihood estimation}
\label{sec.mle}

We will consider now a maximum likelihood approach to the estimation of the state $\rho$. 
Let us recall the terms of the problem: we are given a sequence 
$Y_1,Y_2 \dots, Y_n $ of i.i.d. random variables $Y_i=(X_i,\Phi_i)$ with values in 
$\R\times[0,\pi]$ and probability density $p_\rho$ depending on the parameter $\rho$ which 
is an infinite dimensional density matrix. We would like to find an estimator of $\rho$
\begin{equation*}
\bar{\rho}^{(n)}=\bar{\rho}^{(n)}(Y_1,Y_2 ,\dots, Y_n).
\end{equation*} 
Let $\rho$ and $\tau$ be density matrices and denote by
\begin{equation}\label{eq.hellinger.distance}
h(P_\rho,P_\tau):=\left(\int(\sqrt{p_\rho}-\sqrt{p_\tau})^2~dxd\varphi\right)^{1/2},
\end{equation}
the Hellinger distance between the two probability distributions. 
The following relations are well known
\begin{equation}\label{eq.inequalities.hellinger.norm_1}
d_{\text{tv}}(P_\rho,P_\tau)=\frac{1}{2}\|p_\rho-p_\tau\|_1\leq h(P_\rho,P_\tau)\leq 
\sqrt{\|p_\rho-p_\tau\|_1},
\end{equation}
An important property which is true for any measurement is the following inequality between the 
classical and quantum distances, cf.~\cite{Guta},
\begin{equation}\label{eq.norm_1.inequality}
\|p_\rho-p_\tau\|_1\leq \|\rho-\tau\|_1.
\end{equation}
From (\ref{eq.inequalities.hellinger.norm_1}) and (\ref{eq.norm_1.inequality}) we obtain
\begin{equation}\label{contractivity.T}
h(P_{\rho}, P_{\tau})\leq \sqrt{\|\rho-\tau\|_1},
\end{equation}
for arbitrary states $\rho,\tau $. As we have seen previously, the inverse map from $P_\rho$ to $\rho$ is unbounded thus we do not have an inequality in the opposite direction to the one above. However we can prove the continuity of the inverse map by using a matrix analogue of the Scheff\'{e}'s lemma from classical probability (see \cite{Williams}) stating that if a sequence of probability densities converges pointwise almost everywhere to a probability density, then they also converge in total variation norm. The matrix Scheff\'{e}'s lemma which can be found in \cite{Simon.3} says that 
if $\rho, \rho^{(1)}, \rho^{(2)},\dots$ are density matrices (positive and trace one) and if the coefficients $\rho^{(n)}_{i,j}$ converge to $\rho_{i,j}$ as $n$ goes to infinity, for any fixed indices $i,j$, then $\| \rho^{(n)}-\rho\|_1 \to 0$. 
But by equation (\ref{eq.quantum_tomographic_rho_n,n+d}) if the sequence 
$P_{\rho^{(n)}}$ converges to $P_\rho$ as $n\to\infty$ with respect to the 
$d_{\text{tv}}$-distance then $\rho^{(n)}_{i,j}$ converges to $\rho_{i,j}$ and thus 
$\| \rho^{(n)}-\rho\|_1 \to 0$, completing the proof of the continuity of the map from 
$P_\rho$ to $\rho$. In particular we have $\| \rho^{(n)}-\rho\|_2 \to 0$ due to the inequality 
between the $L_1$ and $L_2$ norms
\begin{equation}
\|\tau-\rho\|_2 \leq \|\tau-\rho\|_1.
\end{equation}

The maximum likelihood estimator of $\rho$ is defined as 
\begin{equation}
\text{arg}\max_{\tau}
\sum_{\ell=1}^n\log p_\tau(X_\ell,\Phi_\ell)
\end{equation}
where the maximum is taken over all density matrices in on the space $L_2(\R)$. 
However there exist density matrices $\tau$ such that the probability density $p_\tau$ 
takes arbitrarily high values at all the points $(X_\ell,\Phi_\ell)$. To see this let us first 
remind the reader that any density matrix is a convex combination of ``pure states'' which are 
projections $\mathbf{P}(\psi)$ on one dimensional sub-spaces of $L_2(\R)$ generated by vectors $\psi$ 
which can be written as a Fourier sum 
\begin{equation}
\psi(x)=\sum_{k=0}^\infty\alpha_k\psi_k(x)
\end{equation}
in the basis $\{\psi_k\}$ given in equation (\ref{eq.psi_n}), 
with $\int|\psi(x)|^2dx=\sum_k|\alpha_k|^2=1$. 
For any such state the corresponding probability density is 
\begin{equation}\label{eq.probability.pure.state}
p_\psi(x,\phi)=\Big|\sum_{k=0}^{\infty}\alpha_k e^{ik\phi}\psi_k(x)\Big|^2= 
|\psi^{\phi}|^2(x),
\end{equation}
where $\psi^{\phi}(x)$ is the square integrable function with Fourier 
coefficients $e^{ik\phi}\alpha_k$. It is clear that there exists a one-to-one relation between 
$\psi$ and $\psi^{\phi}$ which preserves the $L_2$-norms, thus we can choose vectors 
$\varphi_1,\dots \varphi_n$ such that $\int|\varphi_\ell(x)|^2dx=1$ and 
$|\varphi_\ell^{\Phi_\ell}(X_\ell)|^2>C$ for all $\ell=1,\dots n$ and arbitrary $C>0$. 
Then the density matrix  
\begin{equation}
\rho=\frac{1}{n}\sum_{k=1}^n \mathbf{P}(\varphi_k) 
\end{equation}
representing a statistical mixture of the pure states leads to the likelihood 
\begin{equation}
p_{\rho}(x,\phi)= \prod_{\ell}^n
\Big( \frac{1}{n}\sum_{k=1}^n |\varphi_k^{\Phi_\ell} (X_\ell)|^2 \Big)
\geq \Big(\frac{C}{n}\Big)^n.
\end{equation}
which can be arbitrarily high for any fixed $n$. 
This drawback can be corrected by using for example penalized maximum likelihood estimators or by restricting the state space to some subspace $\mathcal{Q}(n)$ of density matrices such that for any amount of data the maximum of the likelihood over $\mathcal{Q}(n)$ exists and $\cup_{n}\mathcal{Q}(n)$ is dense in the space of all density matrices $\mathcal{S}$ with respect to some chosen distance function. 
Such a method is called sieve maximum likelihood and we refer to \cite{Geer} and \cite{Wong&Shen} for the general theory. 
The choice of the sieves $\mathcal{Q}(n)$ should be tailored according to the problem one wants to solve, the class of states one is interested in, etc. Here  we will use the same subspaces as for the projection estimator of the previous subsection, that is $\mathcal{Q}(n)$ consists of those states with maximally $N(n)-1$ photons described by density matrices over the subspace spanned by the basis vectors $\psi_0,\dots,\psi_{N(n)-1}$. We will call $\{\mathcal{Q}(n)\}_{n\geq 0}$ the {\it number states sieves} and the dimension $N(n)$ will be an increasing function of $n$ which will be fixed later so as to guarantee consistency.  
\begin{equation}\label{def.Q(n)}
\mathcal{Q}(n)=\left\{\tau\in\mathcal{S}~:~ \tau_{j,k}=0 ~\text{for all}~ 
j\geq N(n) ~\text{or}~ k\geq N(n)\right\}. 
\end{equation}
Notice that the dimension of the space $\mathcal{Q}(n)$ is $N(n)^2-1$. Let now the estimator be
\begin{equation}\label{def.SMLE}
\bar{\rho}^{(n)}:=\text{arg}\max_{\tau\in\mathcal{Q}(n)}
\sum_{\ell=1}^n\log p_\tau(X_\ell,\Phi_\ell),
\end{equation}
where the maximum can be shown to exist for example by using compactness arguments. We will denote the corresponding sieve in the space of probability densities by 
 \begin{equation}
\mathcal{P}(n)=\left\{p_\tau~:~ \tau\in\mathcal{Q}(n)\right\}. 
\end{equation}

The theory of M-estimators \cite{Geer} tells us that the consistency of ML estimators depends essentially on the ``size'' of the parameter space, in our case the sieves $\mathcal{Q}(n)$ or $\mathcal{P}(n)$, which is measured by entropy numbers with respect to some distance, for example the $L_1$-norm on density matrices or the Hellinger distance between probability distributions.  
\begin{definition}
Let $\mathcal{G}$ be a class of density matrices. Let $N_{B,1}(\delta,\mathcal{G})$ be the smallest value of $p\in\N$ for which there exist pairs of Hermitian matrices (not necessarily density matrices) $[\tau_j^L,\tau_j^U]$ with $j=1,\dots,p$ such that $\|\tau_j^L-\tau_j^U\|_1\leq\delta$ for 
all $j$, and such that for each $\tau \in\mathcal{G}$ there is a $j=j(\tau)\in\{1,\dots, p\}$ 
satisfying
\begin{equation*}
\tau_j^L\leq \tau \leq \tau_j^U.
\end{equation*}
Then $H_{B,1}(\delta, \mathcal{G})=\log N_{B,1}(\delta, \mathcal{G})$ is called 
{\it $\delta$-entropy with bracketing} of $\mathcal{G}$.
\end{definition} 
We note that this definition relies on the concept of positivity of matrices and the existence of the $L_{1}$-distance between states. But the same notions exist for the space of integrable functions thus by replacing density matrices with probability densities and selfadjoint operators with functions we obtain the definition of the $\delta$-entropy with bracketing $H_{B,1}(\delta, \mathcal{F})$ for some space of probability densities $\mathcal{F}$, see \cite{Geer}. 
Moreover by using equation (\ref{eq.norm_1.inequality}) and the fact that the linear extension of the map from density matrices to probability densities sends a positive matrix to a positive function, we get that for any $\delta$-bracketing $[\tau_j^L,\tau_j^U]$ for $\mathcal{Q}(n)$, the corresponding functions $[p_j^L,p_j^U]$ form a $\delta$-bracketing for $\mathcal{P}(n)$, i.e. they satisfy $\|p_j^L-p_j^U\|_1\leq \delta$ and for any $p\in\mathcal{P}(n)$ there exists a $j=j(p)$ such that $p_j^L\leq p \leq p_j^U$. Thus 
\begin{equation*}
H_{B,1}(\delta, \mathcal{P}(n))\leq H_{B,1}(\delta, \mathcal{Q}(n)).
\end{equation*}
The following proposition gives an upper bound of the ``quantum'' bracketing 
entropy and in consequence for $H_{B,1}(\delta, \mathcal{P}(n))$. 
Its proof can be found in \cite{Guta} and relies on choosing a maximal number of 
nonintersecting balls centered in $\mathcal{Q}(n)$ having radius $\frac{\delta}{2N(n)}$ and then providing a pair of brackets for each ball. 
\begin{proposition}
Let $\mathcal{Q}(n)$ be the class of density matrices of dimension $N(n)$. Then
\begin{equation}\label{eq.bound.bracheting.entropy.norm1}
H_{B,1}(\delta,\mathcal{Q}(n))\leq C N(n)^2\log\frac{N(n)}{\delta}~. 
\end{equation}
for some constant $C$ independent of $n$ and $\delta$.
\end{proposition}
By combining the previous inequalities with equation (\ref{contractivity.T}) 
we get the following bound for the bracketing entropy of the class of square-root densities 
\begin{equation}
\mathcal{P}^{1/2}(n)=\Big\{ \sqrt{p}_\rho~:~ p_\rho\in \mathcal{P}(n)\Big\},
\end{equation}  
with respect to the $L_2$-distance
\begin{equation}\label{eq.bracketing.entropy.Hellinger}
H_{B,2}(\delta,\mathcal{P}^{1/2}(n))\leq  C_2 N(n)^2\log\frac{N(n)}{\delta}~.
\end{equation}
We will concentrate now on the Hellinger consistency of the sieve maximum likelihood 
estimator $\bar{P}_n$. We will appeal to a theorem from \cite{Geer}, which is similar to 
other results in the literature on non-parametric $M$-estimation 
(see for example \cite{Wong&Shen}). 
There are two competing factors which contribute to the convergence of $h(\bar{P}_n,P_\rho)$. The first is related to the approximation properties of the sieves with respect to the whole parameter space. Such a distance from $\rho$ to the sieve $\mathcal{Q}(n)$ can take different expressions, for example in terms of the 
$\chi^2$-distance between the corresponding probability measures
\begin{equation}
\tau_n:=\underset{p_n\in\mathcal{P}_n}{\mathrm{argmin}}~\chi^2(p_\rho,p_n),
\end{equation}
where the $\chi^2$-distance between two probability distributions is given by
\begin{equation}
\chi^2(P_1,P_2):=\begin{cases}
\int (\frac{dP_1}{dP_2}-1)^2 dP_2 & P_1\ll P_2,\\
\infty& \mathrm{otherwise}. 
\end{cases}
\end{equation}

Notice that $\tau_n$ depends on $n$ through the growth rate of the sieve $N(n)$. The second factor influencing the convergence of $h(\bar{P}_n,P_\rho)$ is the size of the sieves which is expressed by the bracketing entropy. The non-parametric sieve maximum likelihood estimation theory shows that consistency holds if there exists a sequence $\delta_n\to 0$ such that the following {\it entropy integral inequalities} are satisfied for all $n$
\begin{equation}\label{eq.entropy.integral.inequality}
J_B(\delta_n, \mathcal{P}^{1/2}(n)):=\int_{\delta_n^2/c}^{\delta_n} ~
H_B^{1/2}(u, \mathcal{P}(n)^{1/2})du~\vee~\delta_n \leq \sqrt{n}\delta_n^2/c.
\end{equation} 
where $c$ is some universal constant, \cite{Geer}.  
From (\ref{eq.bracketing.entropy.Hellinger}) we get
\begin{equation}\label{eq.entropy.integral}
J_{B}(\delta_n, \mathcal{P}^{1/2}(n))=
O\left[N(n)\delta_n\left(\log \frac{N(n)}{\delta_n}\right)^{1/2}\right],
\end{equation}
which implies the following constraint for $N(n)\to\infty$ and $\delta_n\to 0$ ,
\begin{equation}\label{eq.rate.delta_n}
\frac{N(n)}{\delta_n}=O\left(\sqrt{\frac{n}{\log n}}\right).
\end{equation}
\begin{theorem}\label{th.Hellinger.consistency}
Suppose that the state $\rho$ satisfies $\tau_n\to 0$. 
Let $\bar{\rho}^{(n)}$ be the sieve MLE with $N(n)$ and $\delta_n$ satisfying (\ref{eq.entropy.integral.inequality}), then
\begin{equation}
\mathbb{P}(h(\bar{P}_n,P_\rho)\geq \delta_n+\tau_n)\leq c~ \mathrm{exp}
(-n\delta_n^2/c^2) .
\end{equation}
\end{theorem}

\noindent\textit{Proof.} 
Details can be found in \cite{Guta} based on Theorem 10.13 of \cite{Geer}. 

\qed

From the physical point of view, we are interested in the convergence of the state estimator $\bar{\rho}^{(n)}$ with respect to the $L_1$ and $L_2$-norms on the space of density matrices. Clearly the rates of convergence for such estimators are slower than those of their corresponding probability densities. As shown in the beginning of this subsection the map sending probability densities $p_\rho$ to density matrices $\rho$ 
is continuous, thus an estimator $\bar{\rho}_n$ taking values in the space of density matrices $\mathcal{S}$ is consistent in the $L_1$ or $L_2$-norms if and only if $\bar{P}_n$ converges to $P_\rho$ almost surely with respect to the Hellinger distance. 
\begin{corollary}\label{cor.consistency.mle.densitymatrix}
The Hellinger consistency of $\bar{P}_n$ is equivalent to the $\|\cdot\|_1$-consistency of 
$\bar{\rho}^{(n)}$. In particular, if $\sum_n \mathrm{exp}(-n\delta_n^2/c^2)<\infty$ and the assumptions of Theorem \ref{th.Hellinger.consistency} hold, then we have $\|\bar{\rho}^{(n)}-\rho\|_1\to 0$, a.s..  
\end{corollary}

\subsection{Wigner function estimation} 
\label{sec.Wigner}

The Wigner function plays an important role in quantum optics as an alternative way of 
representing quantum states and calculating an observable's expectation: for any observable $\mathbf{X}$ there exists a function $W_\mathbf{X}$ from $\R^2$ to $\R$ such that 
\begin{equation}
\mathrm{Tr}(\mathbf{X}\rho)=\iint W_\mathbf{X}(q,p)W_\rho(q,p)dqdp.
\end{equation}
Besides, physicists are interested in estimating the Wigner function for the purpose of identifying features which can be easier visualized than read off from the density matrix, for example a ``non-classic'' state may be recognized by its patches of negative Wigner function, while ``squeezing'' is manifest through the oval shape of the support of the Wigner function, see Table \ref{tbl.states} and Figure \ref{fig.examples}. 
As described in Subsection \ref{sec.qhomodyne} the Wigner function should be seen formally as a joint density of the observables $\mathbf{Q}$ and $\mathbf{P}$ which may take non-negative values, reflecting the fact that the two observables cannot be measured simultaneously. However the Wigner function shares some common properties with probability densities, in particular their marginals $\int W_{\rho}(q,p) dq$ and $\int W_{\rho}(q,p) dp $ are probability densities on the line. In fact this is true for the marginals in any direction which are nothing else then the densities $p_\rho(x,\phi)$. On the other hand there exist probability densities which are not Wigner functions and vice-versa, for example the latter cannot be too ``peaked'':
\begin{equation}\label{eq.uniform.bounded.Wigner}
| W_{\rho}(q,p)|  \leq \frac{1}{\pi}, \qquad \mathrm{for~all}~(q,p)\in\R^2,~ 
\rho\in\mathcal{S}.
\end{equation}
As a corollary of this uniform boundedness we get
\begin{equation}\label{eq.norm1.norminfty.inequality}
\| W_{\rho}-W_{\tau}\|_\infty\leq \frac{1}{\pi}\|\rho-\tau\|_1.
\end{equation} 
for any density matrices $\rho$ and $\tau$. Indeed we can write $\rho-\tau=\rho_+-\rho_-$ where 
$\rho_+$ and $-\rho_-$ represent the positive and negative part of $\rho-\tau$. Then  
\begin{eqnarray}
&&\| W_{\rho}-W_{\tau}\|_\infty=\| W_{\rho_+}-W_{\rho_-}\|_\infty\leq 
\| W_{\rho_+}\|+\|W_{\rho_-}\|_\infty\leq\nonumber\\
&& \frac{1}{\pi}(\|\rho_+\|_1+\|\rho_-\|_1)=\frac{1}{\pi}\|\rho-\tau \|_1\nonumber.
\end{eqnarray}
Another important property is the fact that the linear span of the Wigner functions is dense in $L_2(\R^2)$, the 
space of real valued, square integrable functions on the plane, and there is an isometry 
(up to a constant) between the space of Wigner functions and that of density 
matrices with respect to the $L_2$-distances
\begin{equation}\label{eq.isometry}
\| W_{\rho}- W_{\tau}\|_2^2 =:\iint |W_{\rho}(q,p)- W_{\tau}(q,p)|^2dp dq=
\frac{1}{2\pi}\parallel \rho-\tau \parallel_2^2. 
\end{equation} 
In Section \ref{sec.problem} we have described the standard estimation method employed in computerized tomography which used a regularized kernel $K_c$ with bandwidth $h_n=1/c$ converging to zero as $n\to\infty$ at an appropriate rate. This type of estimators for the Wigner function will be analyzed separately in future work in the minimax framework along the lines of \cite{Cavalier}. The estimators which we propose in this subsection are of a different type, they are based on estimators for $\rho$ plugged into the following linearity equation
$$ 
W_{\rho}(q,p) = \sum_{k,j} \rho_{k,j} W_{k,j}(q,p),
$$
where $W_{k,j}$ are known functions corresponding to the matrix with the 
entry $(k,j)$ equal to $1$ and all the rest equal to  zero, see \cite{Leonhardt}. 
The isometry (\ref{eq.isometry}) implies that the family $\{W_{k,j}\}_{k,j=0}^{\infty}$ forms an orthogonal basis of $L_2(\mathbb{R}^2)$. 
Following the same idea as in the previous section we consider the projection estimator
$$
\widehat{W}^{(n)}(q,p) = \sum_{k,j=0}^{N(n)-1} \left( \frac{1}{n} 
\sum_{\ell=1}^{n} F_{k,j}(X_\ell,\Phi_\ell) \right) W_{k,j}(q,p).
$$
\begin{corollary}\label{th.wigner.norm2} Let $N(n)$ be such that $N(n)\to\infty$ and 
\begin{equation*}
N(n)=o(n^{3/7}). 
\end{equation*}
Then
\begin{equation*}
\mathbb{E}\|\widehat{W}^{(n)}-W_\rho\|_2^2=\|W_{\rho}^{(n)}-W_\rho\|_2^2+
 O\left(\frac{N(n)^{7/3}}{n}\right).
\end{equation*}
\end{corollary}

\noindent{\it Proof:} Apply isometry property and Theorem \ref{th.norm2}.

\qed

Similarly we can extend the SML estimator of the density matrix to the Wigner function. 
Define the subspace
$$ 
\mathcal{W}(n) = \{ W_{\rho}: \, \rho \in \mathcal{Q}(n)\},
$$
with $\mathcal{Q}(n)$ as in equation (\ref{def.Q(n)}), and define the corresponding SML 
estimator as $\widebar{W}^{(n)} = W_{\bar{\rho}^{(n)}}$ where $\bar{\rho}^{(n)}$ was defined in 
(\ref{def.SMLE}).
\begin{corollary}\label{th.Hellinger.consistency.wigner}
Suppose that $\rho$ satisfies $\tau_n\to 0$. Let $\widebar{W}_{\rho}^{(n)}$ be the SML estimator with $N(n)$ and $\delta_n$ satisfying (\ref{eq.entropy.integral.inequality}) 
and $\sum_n \mathrm{exp}(-n\delta_n^2/c^2)<\infty$. Then we have 
$$\|\widebar{W}^{(n)}-W_\rho\|_2\to 0$$ almost surely. Under the same conditions
$$\|\widebar{W}^{(n)}-W_\rho\|_\infty\to 0$$ almost surely.  
\end{corollary}

\noindent{\it Proof:} Apply the inequalities (\ref{eq.norm1.norminfty.inequality}, 
\ref{eq.isometry}) and Corollary \ref{cor.consistency.mle.densitymatrix}. 

\qed

\section{Noisy observations}\label{sec.efficiency}

The homodyne tomography measurement as presented up to now does not take into account various losses (mode mismatching, failure of detectors) in the detection process which modify the distribution of results in a real measurement compared with the idealized case.  Fortunately, an analysis of such losses (see \cite{Leonhardt}) shows that they can be quantified by a single {\it efficiency} coefficient $0<\eta<1$ and the change in the observations amounts replacing $X_i$ by the noisy observations
\begin{equation}
X_i':=\sqrt{\eta}X_i+\sqrt{(1-\eta)/2}\xi_i
\end{equation}  
with $\xi_i$ a sequence of i.i.d.~standard Gaussians which are independent of all $X_j$. 
The problem is again to estimate the parameter $\rho$ from $Y_i'=(X_i',\Phi_i),$ for 
$i=1,\dots,n$. The efficiency-corrected probability density is then the convolution 
\begin{equation}
p_\rho(y,\phi;\eta)=(\pi(1-\eta))^{-1/2}
\int_{-\infty}^\infty p(x,\phi)\mathrm{exp}\left[-\frac{\eta}{1-\eta}
(x-\eta^{-1/2}y)^2\right] ~dx.
\end{equation}

The physics of the detection process detailed in \cite{Leonhardt} offers an alternative route 
from the state to the probability density of the observations $Y_i'$. 
In a first step one performs a {\it Bernoulli transformation} $B(\eta)$ on the state 
$\rho$ which is a quantum equivalent of the convolution with noise for probability densities, and obtains a new density matrix $\rho^\eta$. 
To understand the Bernoulli transformation let us consider first the diagonal elements 
$\{ p_k=\rho_{k,k},~k=0,1..\}$ and $\{q_j=\rho^{\eta}_{j,j},~j=0,1..\}$ which are 
both probability distributions over $\N$ and represent the statistics of the number of 
photons in the two states. Let $ b_k^{k+p}=\binom{k+p}{k}\eta^k(1-\eta)^p$ be the 
binomial distribution. Then
\begin{equation}
q_j=\sum_{k=j}^\infty b_j^{k}(\eta)p_k
\end{equation}  
which has a simple interpretation in terms of an ``absorption'' process by which each 
photon of the state $\rho$ goes independently through a filter and is allowed to pass with probability $\eta$ or absorbed with probability $1-\eta$. The formula of the Bernoulli transformation for the whole matrix is 
\begin{equation}\label{eq.Bernoulli}
\rho^{\eta}_{j,k}=
\sum_{p=0}^{\infty}\left[ b_j^{j+p}(\eta) b_k^{k+p}(\eta)\right]^{1/2}\rho_{j+p,k+p}.
\end{equation}  
The second step is to perform the usual quantum tomography measurement with ideal detectors on the ``noisy'' state $\rho^\eta$ obtaining the results $Y_i'$ with density $p_\rho(x,\phi;\eta)$. It is noteworty that the transformations $B(\eta)$ form a semigroup, that is they can be composed as $B(\eta_1)B(\eta_2)=B(\eta_1\eta_2)$ and the inverse of $B(\eta)$ is simply obtained by replacing $\eta$ with $\eta^{-1}$ in equation (\ref{eq.Bernoulli}). Notice however that if $\eta\leq 1/2$ the power series $(1-\eta^{-1})^k$ appearing in the inverse transformation 
diverges, thus we need to take special care in this range of parameters.

A third way to compute the inverse map from $p_\rho(x,\phi;\eta)$ to $\rho$ is by using 
pattern functions depending on $\eta$ which incorporate the deconvolution 
map from $p_\rho(x,\phi;\eta)$ to $p_\rho(x,\phi;1)$: 
\begin{equation}\label{eq.quantum_tomographic_rho_n,n+d_eta}
\rho_{k,j}=\int_{-\infty}^\infty dx \int_0^\pi \frac{d\phi}{\pi} p_\rho(x,\phi;\eta)
f_{k,j}(x;\eta). 
\end{equation}    
Such functions are analyzed in \cite{D'Ariano.1,D'Ariano.2} where it is argued that the method has a fundamental limitation for $\eta\leq 1/2$ in which case the 
pattern functions are unbounded, while for $\eta>1/2$ numerical calculations show 
that their range grows exponentially fast with both indices $j,k$. 

The two estimation methods considered in Section \ref{sec.DensityMatrixEstimation} can be 
applied to the state estimation with noisy observations. The projection estimator has the same 
form as in Subsection \ref{sec.pfe} with a similar analysis of the mean $L_2$-risk taking 
into account the norms of the new pattern functions $f_{k,j}(x;\eta)$. 
The sieve maximum likelihood estimator follows the definition in Subsection 
\ref{sec.mle} and a consistency result can be formulated on the lines of Corollary 
\ref{cor.consistency.mle.densitymatrix}. We expect however that the rates of convergence will 
be dramatically slower and we will leave this analysis for a separate work.

\section{Experimental results}
\label{sec.experimental}

\begin{figure}[Htbp] 
\begin{center}
\vspace{1cm}
\includegraphics[height=6cm,width=14cm]{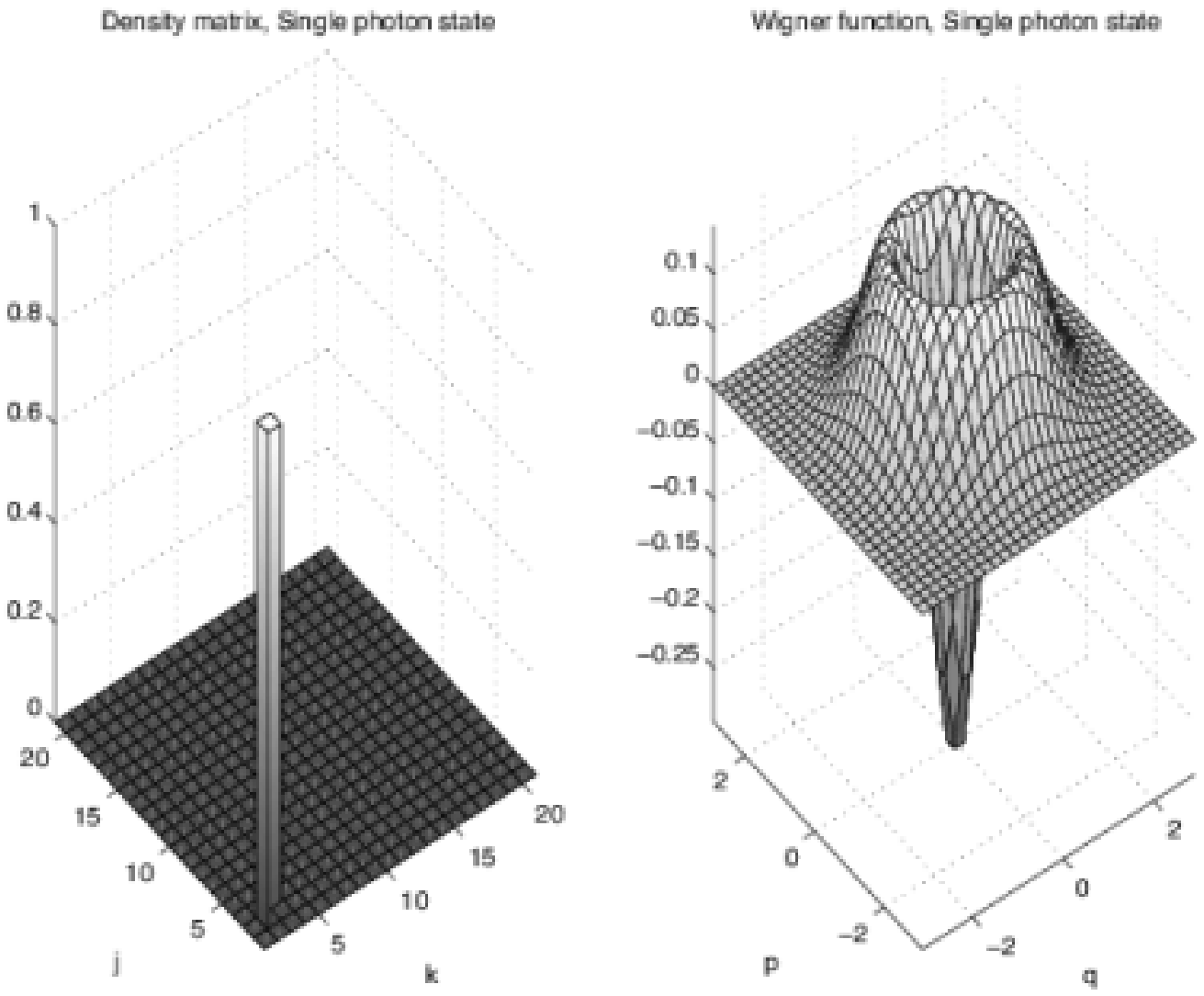}\\[2em]
\includegraphics[height=6cm,width=14cm]{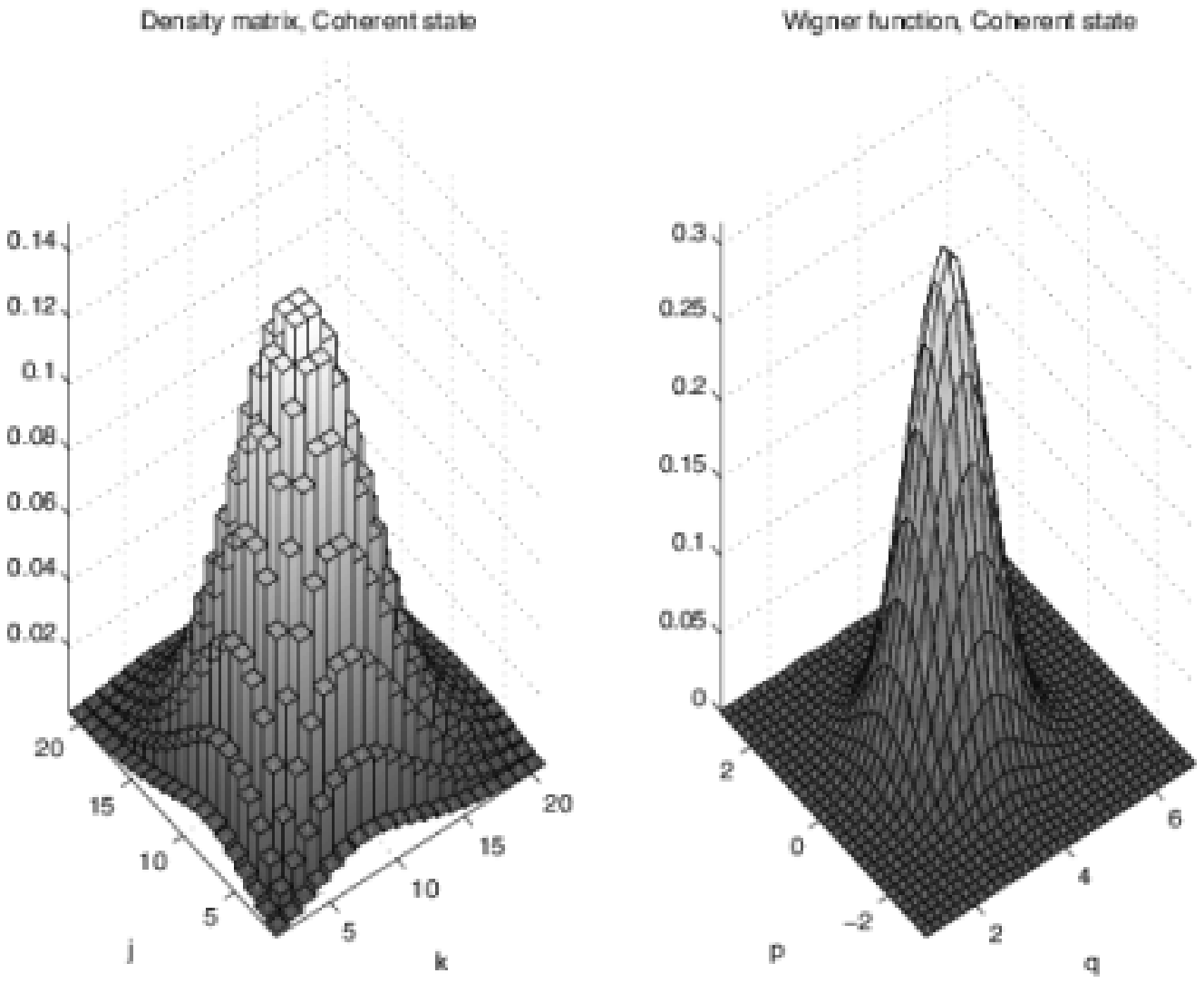}\\[2em]
\includegraphics[height=6cm,width=14cm]{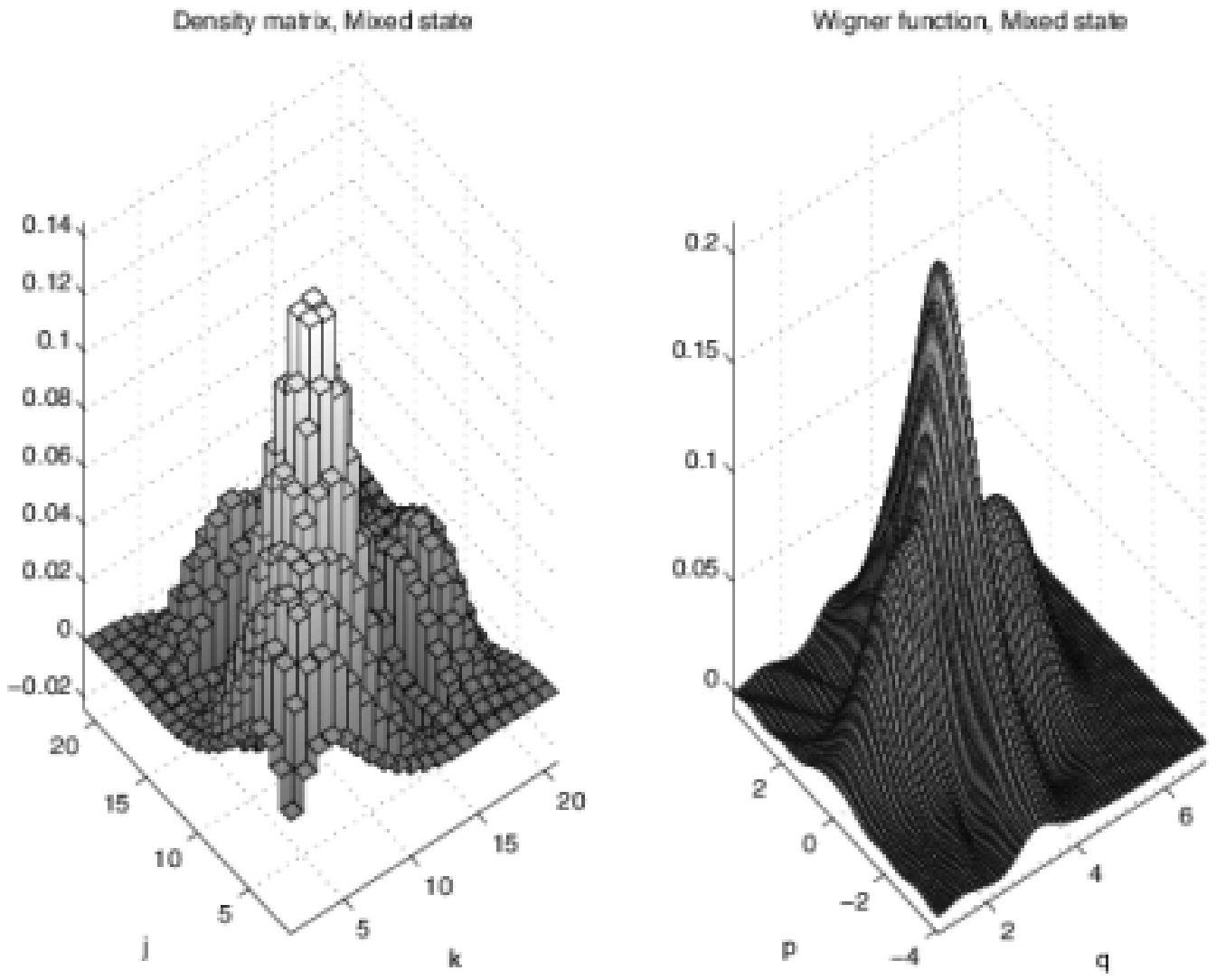}\\
\end{center}
\caption{Graphical representation of quantum states} \label{fig.examples}

\end{figure}

In this section we study the performance of the Pattern Function Projection estimator 
and the Sieve Maximum Likelihood estimator using simulated data. 
In Table \ref{tbl.states} we showed some examples of density matrices and Wigner functions of quantum states. In Figure \ref{fig.examples}, we display their corresponding graphical representation. For some of them the corresponding probability distribution can be expressed explicitly and it is possible to simulate data. In particular we shall simulate data from QHT measurements on squeezed states with efficiency $\eta=1$.

\subsection{Implementation}

In order to implement the two estimators we need to compute the basis functions $\psi_n$ and the functions $\phi_n$, which are solutions of Schr\"odinger equation (\ref{eq.Schrodinger}). For this, we use an appropriate set of recurrent equations, see \cite{Leonhardt}, Ch. 5. Pattern functions can then be calculated as
\begin{equation*}
f_{k,j}(x) = 2 x \psi_{j}(x)\phi_{k}(x) - \sqrt{2(j+1)}\psi_{j+1}(x)\phi_{k}(x) - 
\sqrt{2(k+1)}\psi_{j}(x)\phi_{k+1}(x),
\end{equation*}
for all $j \geq k$,  otherwise $f_{k,j}(x) = f_{j,k}(x)$ and then 
$F_{k,j}(x,\theta) = f_{k,j}(x) e^{i(k-j)\theta}$.

On the practical side, finding the maximum of the likelihood function over a set 
of density matrices is a more complicated problem due to the positivity and trace one 
constraints which must be taken into account. 
A solution was proposed in \cite{D'Ariano.6}, where the restriction on positivity of a density matrix $\rho$ is satisfied by writing the Cholevski decomposition 
\begin{equation}
\rho = T^*T
\end{equation}
where $T$ of {\it upper triangular} matrices of the same dimension as $\rho$ with complex coefficients above the diagonal and reals on the diagonal. The normalization condition $\mathrm{Tr}\rho=1$ translates into $\|T\|_2=1$ which 
defines a ball in the space of upper triangular matrices with the $L_2$-distance. 
We will denote by $\mathcal{T}(n)$ the set of such matrices having dimension $N(n)$. 
The sieve maximum likelihood estimator is the solution of the following  optimization problem with $N=N(n)$
\begin{eqnarray}
\widehat{T} & = & \underset{T\in\mathcal{T}(n)}{\arg \max}~~ \sum_{\ell=1}^{n} \log p_{\rho}(X_\ell,\Phi_\ell) 
\nonumber \\
& = & \underset{T\mathcal{T}(n)}{\arg \max}~~\sum_{\ell=1}^{n} \log 
\sum_{k=1}^{N-1} 
\Big|\sum_{m=0}^{k} \nonumber T_{k,m}\psi_{m}(X_\ell) e^{i m \Phi_\ell} \Big|^2.
\end{eqnarray}

The numerical optimization was performed using a classical descendent method with constraints. Notice that we have an optimization problem on $N^2-1$ real variables. 
Given the problem of high dimensionality and computational cost we propose an alternative method to the procedure mentioned above. It exploits the mixing 
properties of our model. Any density matrix of dimension $N$  can be written as convex 
combination of $N$ {\it pure} states, i.e.~$\rho = \sum_{r=0}^{N-1} p_r \rho_r$, where $ p_r\geq0$, $\sum_rp_r=1$ and $\rho_r$ is a one dimensional projection whose 
Cholevski decomposition is of the form $\rho_r = t_r^{*}t_r$ where $t_r$ is the row vector of dimension $N$ on which $\rho_r$ projects, and $t_r^{*}$ is the column vector of the complex conjugate of $t_r$. It should be noted that even though decomposition of our state in pure states is not unique this is not a problem given we are actually not interested in this representation but in the resulting convex combination, the state itself. Now we can state the problem as to find the maximizer of the loglikelihood
\begin{equation*}
L(\{(X_{\ell},\Phi_{\ell})\};p,t )= \sum_{\ell=1}^{n} 
\log \sum_{r=0}^{N-1} p_r p_{\rho_r}(X_\ell,\Phi_\ell) =  
\sum_{\ell=1}^{n} \log  \sum_{r=0}^{N-1}
 p_r \Big|\sum_{m=0}^{N(n)-1} \nonumber t_{r,m}\psi_{m}(X_\ell) e^{i m \Phi_\ell} \Big|^2,
\end{equation*}
where $t_{r,m}$ represents the m$^{\text{th}}$ coordinate of $t_r$. 
We now maximize over all $\beta\in \mathrm{B}$ where
\begin{equation*}
\mathrm{B}= \left\{ \{p_r, t_r\}_{r=0}^{N-1}~:~ p_r \geq 0, 
\sum_{r=0}^{N-1} p_r =1, \|t_r\|_2=1, \text{ for~all}~ r\right\}.
\end{equation*}

We propose an EM algorithm as an alternative method to the one presented in \cite{D'Ariano.6}. 
See, \cite{Dempster} for an exposition on the formulation of the EM algorithm for problems of 
ML estimation with mixtures of distributions. The iteration procedure is given then in the 
following steps.\\
\noindent 1) Compute the expectation of the conditional likelihood:
$$
Q(\beta|\beta^{\text{old}}) = 
\frac{1}{n} \sum_{\ell=1}^{n} \sum_{r=0}^{N(n)-1} 
\frac{p_r^{\text{old}}p_{\rho_r}^{\text{old}}(X_\ell,\Phi_\ell) }{\sum_{j=0}^{N(n)-1} p_j^{\text{old}}p_{\rho_j}^{\text{old}}(X_\ell,\Phi_\ell) } \log p_r p_{\rho_r}(X_\ell,\Phi_\ell). 
$$
\noindent 2) Maximize $Q(\beta|\beta^{\text{old}})$ over all $\beta\in \mathrm{B}$ and obtain $\beta^{\text{new}}$ with components 
\begin{eqnarray}
&&p_r^{\text{new}} = \frac{1}{n} \sum_{\ell=1}^n \frac{p_r^{\text{old}}p_{\rho_r}^{\text{old}}(X_\ell,\Phi_\ell) }{\sum_{j=0}^{N(n)-1} p_j^{\text{old}}p_{\rho_j}^{\text{old}}(X_\ell,\Phi_\ell) },\nonumber \\
&&\rho_{r}^{\text{new}} = t_r^*t_r, \qquad
t_r  =  \arg \max_{\|t\|_2=1} \sum_{\ell=1}^{n}  f_{r,\ell}^{\text{old}} \log p_{\rho_r}(X_\ell,\Phi_\ell),  \nonumber
\end{eqnarray}
where
$$
f_{r,\ell}^{\text{old}} = \frac{p_r^{\text{old}}p_{\rho_r}^{\text{old}}(X_\ell,\Phi_\ell) }{\sum_{j=0}^{N(n)-1} p_j^{\text{old}}p_{\rho_j}^{\text{old}}(X_\ell,\Phi_\ell) }.
$$

As initial condition one could take very simple ad hoc states. For example, take $t_0$ as the vector $(1,0,\ldots,0)$ and $p_0=1$ and $t_r$ to be the null vector and $p_r=0$ for $r>0$. This corresponds to the {\em one photon} state. 
Another possible combination is to take $t_{r,m}=\delta_r^m$ and $p_r=\frac{1}{N}$. This corresponds to the state represented by a diagonal matrix, called a chaotic state. Another strategy is to consider a preliminary estimator, based on just few observations, diagonalize it and take $t_r$ equal to its eigenstates and $p_r$ the corresponding eigenvalues. In this way one hopes to start the iteration from a state closer to the optimum one. In terms of speed our simulations suggests to use the EM version as dimension grows. Establishing any objective comparison between direct optimization 
and EM algorithm has proven to be difficult given the dependency on initial conditions,  and high dimensionality of the problem. 

\subsection{Analysis of results}

\begin{figure}[Hhbt]
\begin{center}
\includegraphics[width=14cm]{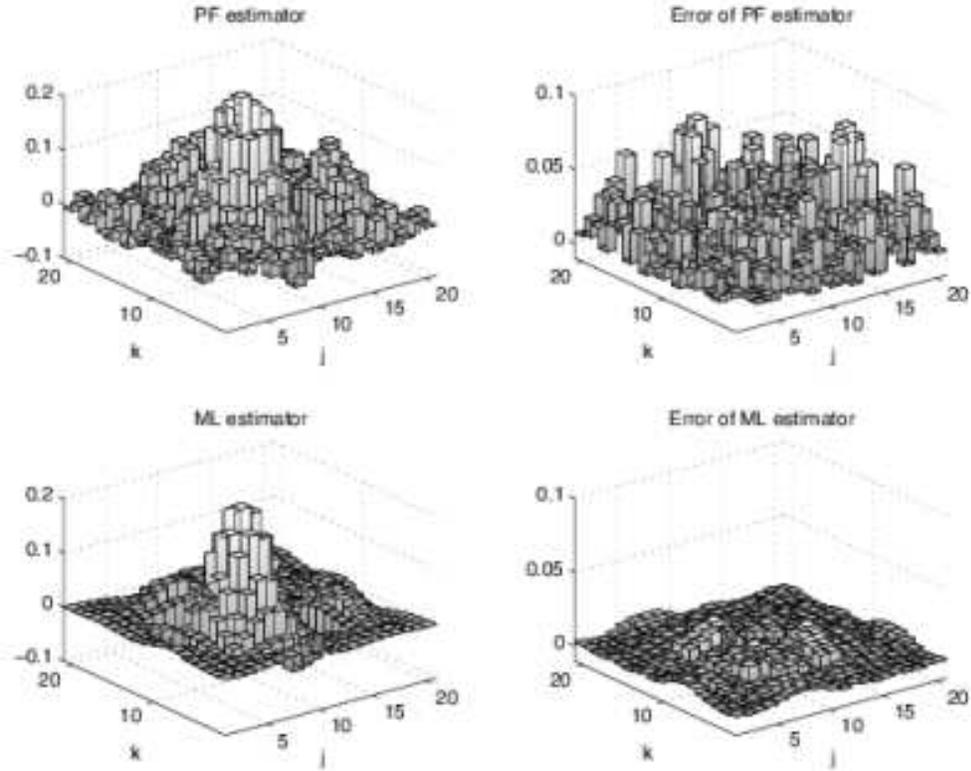}
\caption{Estimation of Squeezed state. First column, from top to bottom, PFP and SML estimators. 
Second column, corresponding errors}\label{fig.estimation.squeezed}
\end{center}
\end{figure}

\begin{figure}[Htb]
\begin{center}
\includegraphics[width=14cm]{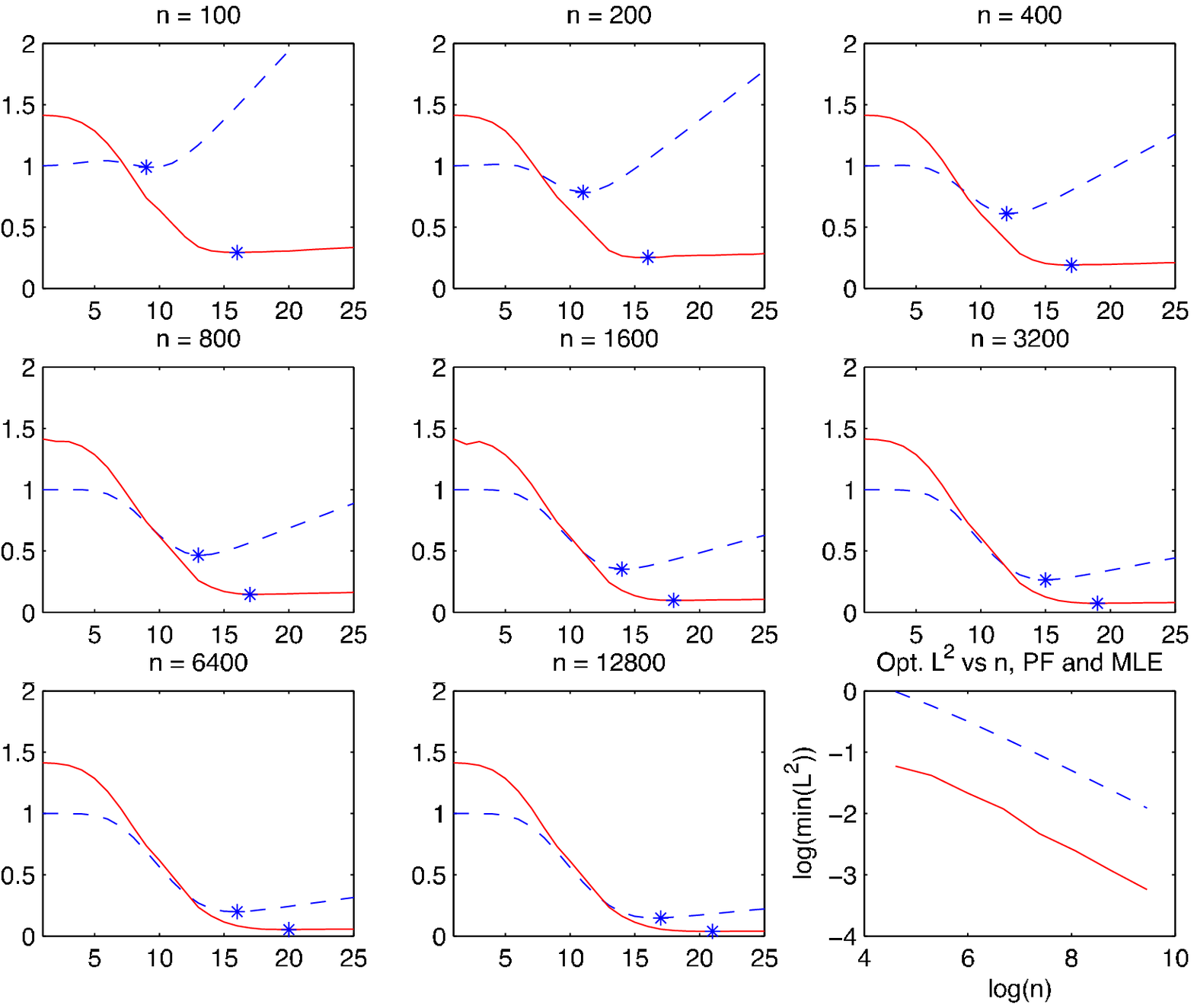}
\caption{$L_2$ error for Pattern Function Estimator and Sieve Maximum Likelihood Estimator 
and different sample sizes: $ n= 100, \ldots, 12800$. Last graphic represents the optimum 
$L_2$ error for different sample size, using a logarithmic scale on both axis.}\label{fig.MLEvsPF.l2error}
\end{center}
\end{figure}

In Figure \ref{fig.estimation.squeezed} we show the result of estimating the squeezed state defined in Table \ref{tbl.states} using samples of size 1600, for both Pattern Function and Maximum Likelihood estimators. At a first glance one can see that the Pattern Function estimator result is rougher when compared to the Maximum Likelihood 
estimator. This is due to the fact discussed in Subsection \ref{sec.pfe} that the 
variance of $F_{k,j}$ increases as a function of $k$ and $j$ as we move away from the diagonal. 
The relation between quality of estimation and dimension of the truncated estimator is seen more clearly in Figure \ref{fig.MLEvsPF.l2error} where the $L_2$-errors of estimating the coherent state is shown for both estimators at different sample sizes. The $*$-dot represents the point of minimum -- and thus, optimum dimension $N^*(n)$ -- for each curve. The curves presented there are the mean $L_2$-error estimated using 15 simulations for each sample size. From there we can see that the optimum PFP estimator for the sample of size $n=1600$ is the one corresponding to $N^*(n)=15$ while the optimum SML would be obtained using the sieve of size $N^*(n)=19$.

\begin{figure}[hbt]
\begin{center}
\includegraphics[width=7cm]{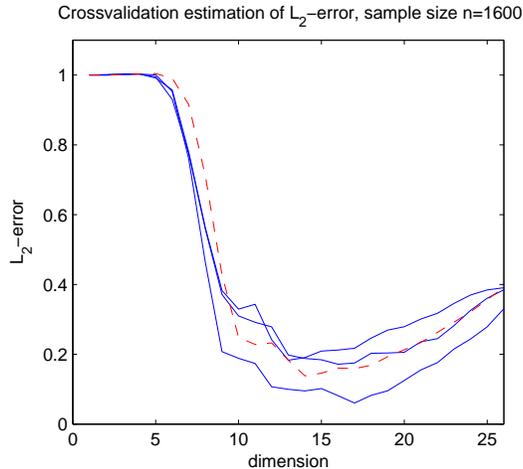}
\caption{Estimated optimal dimension for three simulations with 1600 observations}\label{fig.cross-validation}
\end{center}
\end{figure}

Let us first analyze the performance of PFP estimator. 
Notice that for $N>N^*$ the mean 
$L_2$-error increases quadratically with $N$ due to contribution from the variance term. 
As $n$ increases the variance decreases like $n^{-1}$ for a fixed dimension $N$ and 
consequently, the optimal dimension $N^*(n)$ increases. One can see that the 
minimum is attained rather sharply. This suggests 
that, in order to get a good result a refined method of guessing the optimum 
dimension becomes 
necessary, eg.~BIC, AIC or cross-validation. 
Figure \ref{fig.cross-validation} shows a cross-validation estimator for the $L_2$-error of the PFP estimator for three simulations (continuous lines), each 
one based on 1600 observations with a squeezed state, and for comparison the expected $L_2$-error (dashed line). This represents only a first attempt to implement model selection procedures for this problem which should be investigated more thoroughly from the theoretical and practical point of view.

We pass now to the SML estimator. In Figure \ref{fig.MLEvsPF.l2error} we see that it 
has smaller $L_2$-error the the PFP estimator at its optimum sieve dimension. 
It is remarkable that the behavior of the $L_2$-error, for $N>N^*$ has a different behavior in this case, increasing much slower than the PFP estimator at the right side of its corresponding optimal dimension. This suggests that SML estimators could have a lower risk if the optimum dimension is overestimated. 
The bottom right pane of Figure \ref{fig.MLEvsPF.l2error} shows the optimum value of the $L_2$-risk in terms of sample size. Both axis are represented in a logarithmic scale. The observed linear pattern indicates that the $L_2$-risk decreases as $an^{-\tau}$. The slope of both curves correspond to $\tau \sim 0.4$, showing an almost parametric rate which is not surprising given the smoothness of the example that we consider. The value of $\tau$ for the PFP estimator is a bit smaller than for the SML estimator, confirming its worse performance. Notice also that the constant $a$ is bigger for PFP than for the SML estimator. We expect that the contrast between the two estimators will be 
accentuated when $\eta<1$. 

\begin{figure}[Htbp]
\begin{center}
\includegraphics[height=5cm,width=14cm]{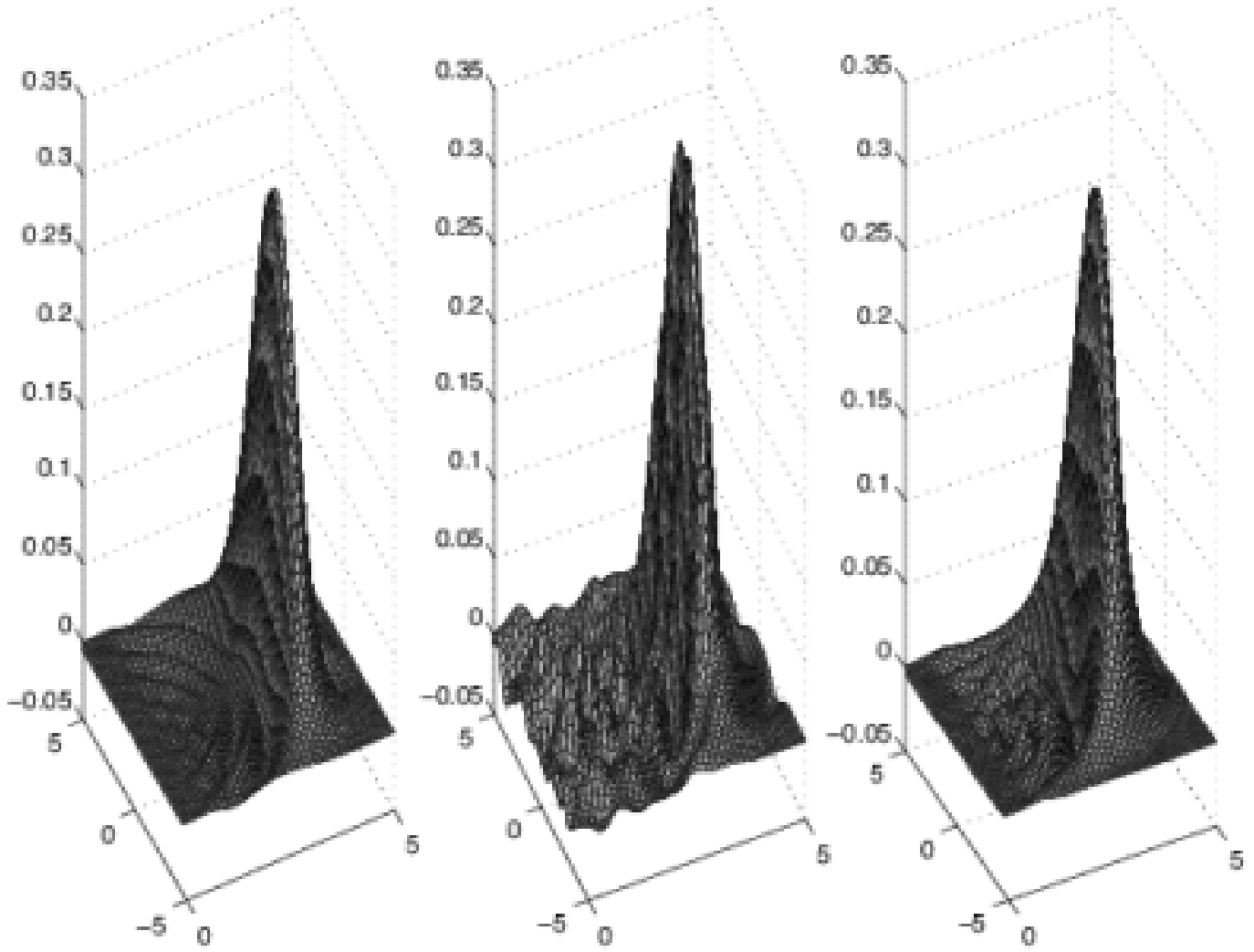}
\end{center}
\caption{Graphical representation of Wigner functions estimation. a) Original Wigner function 
of the squeezed state, b) Pattern Function estimation, c) Maximum Likelihood estimation}
\label{fig.estimation.wigner}
\end{figure}

Finally, in Figure \ref{fig.estimation.wigner} we show the result of estimating the Wigner 
function of the squeezed state using both estimators. As explained in Section 
\ref{sec.Wigner} the corresponding estimator can be obtained 
by plugging-in the density matrix estimator in equation \ref{eq.linearity.Wigner.rho}. 
The density matrix estimators for the same state are represented in Figure 
\ref{fig.estimation.squeezed}.

\section{ Concluding remarks}

In this paper we have proposed a Pattern Function Projection estimator and a 
Sieve Maximum Likelihood estimator for the density matrix of the quantum state and 
its Wigner function. We proved they are consistent for different norms in their corresponding spaces. There are many open statistical questions 
related to quantum tomography and we would like to enumerate a few of them here.

\begin{itemize}

\item {\it Cross-validation.} For both types of estimators, a data 
dependent method is needed for selecting the optimal sieve dimension. 
We mention criteria such as unbiased cross-validation, hard thresholding 
or other types of minimum contrast estimators \cite{Barron&Birge&Massart}.

\item {\it Efficiency $0<\eta<1$.} A realistic detector has detection efficiency $0<\eta<1$ which introduces an additional 
noise in the homodyne data. From the statistical point of view we deal with a 
Gaussian deconvolution problem on top of the usual quantum tomography estimation.

\item {\it Rates of convergence.} Going beyond consistency requires 
the selection of classes of states which are natural both from the physical, as well as 
statistical point of view. One should study optimal and achieved rates of convergence for given classes. For $0<\eta<1$ the rates are expected to be significantly lower than in the ideal case, so it becomes even more crucial to use optimal estimators. In applications, sometimes only the estimation of a functional of $\rho$ such as average number of photons or entropy may be needed. This will require a separate analysis, cf.~\cite{Shen}.

\item{\it Kernel estimators for Wigner function.} When estimating the Wigner function it seems more natural to use a kernel estimator such as in \cite{Cavalier} and to combine this analysis with the deconvolution problem in the case noisy observations $\eta<1$, \cite{Butucea}.

\item{\it Other quantum estimation problems.} 
The methods used here for quantum tomography can be applied in other problems 
of quantum estimation, such as for example the calibration of measurement devices or 
the estimation of transformation of states under the action of quantum mechanical devices.
\end{itemize}

\bibliographystyle{Chicago}

\end{document}